\documentclass[twoside,11pt]{article}
\title{} \author{} \date{}
\usepackage{latexsym,amssymb,times}
\input amssym.def
\usepackage{graphicx}
\newtheorem{te}{Theorem}[section]
\newtheorem{prop}[te]{Proposition}

\newtheorem{fac}[te]{Fact}

\newtheorem{rem}[te]{Remark}
\newtheorem{ex}[te]{Example}

\def\dok{\noindent{\bf Proof. }}
\def\kdok{\hfill $\Box$ \par \vspace*{2mm} }
\def\a{\alpha}
\def\b{\beta}
\def\g{\gamma}

\def\f{\varphi}
\def\p{\psi}
\def\o{\omega}
\def\k{\kappa}
\def\r{\rho}
\def\s{\sigma}

\def\l{\lambda}

\def\P{{\mathbb P}}
\def\Q{{\mathbb Q}}

\def\N{{\mathbb N}}
\def\X{{\mathbb X}}
\def\Y{{\mathbb Y}}
\def\Z{{\mathbb Z}}

\def\A{{\mathbb A}}

\def\BG{{\mathbb G}}
\def\BR{{\mathbb R}}

\def\BL{{\mathbb L}}

\def\CB{{\mathcal B}}

\def\CD{{\mathcal D}}
\def\CS{{\mathcal S}}

\def\CU{{\mathcal U}}

\def\CN{{\mathcal N}}
\def\CM{{\mathcal M}}

\def\CO{{\mathcal O}}
\def\CX{{\mathcal X}}
\def\la{\langle}
\def\ra{\rangle}
\def\ft{{\mathfrak t}}
\def\fm{{\mathfrak m}}
\def\fc{{\mathfrak c}}

\def\dom{\mathop{\mathrm{dom}}\nolimits}
\def\ran{\mathop{\mathrm{ran}}\nolimits}
\def\id{\mathop{\mathrm{id}}\nolimits}

\def\Aut{\mathop{\rm Aut}\nolimits}
\def\Cond{\mathop{\rm Cond}\nolimits}

\def\Bad{\mathop{\rm Bad}\nolimits}
\def\ar{\mathop{\rm ar}\nolimits}

\def\At{\mathop{\mathrm{At}}\nolimits}

\def\Lim{\mathop{\rm {Lim}}\nolimits}

\def\DO{\mathop{\mathrm{DO}}\nolimits}
\def\cov{\mathop{\mathrm{cov}}\nolimits}
\def\add{\mathop{\mathrm{add}}\nolimits}

\def\Int{\mathop{\mathrm{Int}}\nolimits}
\def\Pc{\mathop{\mathrm{Pc}}\nolimits}
\def\Min{\mathop{\mathrm{Min}}\nolimits}

\def\he{\mathop{\mathrm{ht}}\nolimits}
\def\cf{\mathop{\mathrm{cf}}\nolimits}
\def\Is{\mathop{\mathrm{Is}}\nolimits}

\def\Nwd{\mathop{\mathrm{Nwd}}\nolimits}
\def\Borel{\mathop{\mathrm{Borel}}\nolimits}
\def\Fn{\mathop{\mathrm{Fn}}\nolimits}
\def\supp{\mathop{\mathrm{supp}}\nolimits}
\def\cl{\mathop{\mathrm{cl}}\nolimits}
\def\Bad{\mathop{\mathrm{Bad}}\nolimits}
\def\RD{\mathop{\mathrm{D}}\nolimits}
\begin{document}
\thispagestyle{plain}
\begin{center}
           {\large \bf {\uppercase{Back and forth systems witnessing irreversibility }}}
\end{center}
\begin{center}
{\bf Milo\v s S.\ Kurili\'c}\footnote{Department of Mathematics and Informatics, University of Novi Sad,
                                      Trg Dositeja Obradovi\'ca 4, 21000 Novi Sad, Serbia,
                                      e-mail: milos@dmi.uns.ac.rs}
\end{center}
\begin{abstract}
\noindent
If $L$ is a relational language, then an $L$-structure $\X=\la X,\bar \r \ra$ is reversible
iff there is no interpretation $\bar \s\varsubsetneq \bar \r$ such that the structures $\la X,\bar \s \ra$ and $\la X,\bar \r \ra$ are isomorphic.
We show that $\X$ is not reversible
iff there is a back and forth system $\Pi$ of partial self-condensations of $\X$
containing one which is not a partial isomorphism
and having certain closure properties.
Using that characterization we detect several classes of non-reversible partial orders
containing, for example,
homogeneous-universal posets (in particular, the random poset),
the divisibility lattice, $\la \N ,\,\mid\,\ra$,
the ideals $[\k ]^{<\l}$, the meager ideal in the algebra $\Borel (\o ^\o)$,
and the direct powers of rationals, $\Q ^\k$, and integers, $\Z ^\k$.
Some of the results are obtained under additional set-theoretic assumptions.


\end{abstract}
\section{Introduction}\label{S1}
Reversibility  is a phenomenon which was noticed and considered in several categories of mathematical structures and their morphisms.
For example, a topological space $\CX$ is reversible
iff each continuous bijection (condensation) $f:\CX \rightarrow \CX$ is a homeomorphism (i.e., an automorphism of $\CX$)
or, equivalently, an open mapping.
So, a consequence of the Brouwer invariance of domain theorem is that the Euclidean spaces of finite dimension, $\BR ^n$, are reversible,
while by the open mapping theorem (Banach-Schauder theorem) each continuous linear bijection of a Banach space onto itself is a homeomorphism
(thus here we have a restriction to linear mappings; normed spaces of infinite dimension are not reversible, see \cite{Raj}).

An analogous phenomenon can be considered in the class of relational structures.
If $L$ is a relational language,
then an $L$-structure $\X=\la X,\bar \r \ra$ is called reversible
iff each bijective homomorphism (condensation) $f:\X \rightarrow \X$ is an automorphism of $\X$
or, equivalently, iff there is no interpretation $\bar \s\varsubsetneq \bar \r$ such that the structures $\la X,\bar \s \ra$ and $\la X,\bar \r \ra$ are isomorphic.
In this context we consider some natural questions.

First, it is evident that all linear orders are reversible.
Moreover, finite products of well orders are well founded posets without infinite antichains (see \cite{Harz}, p.\ 250)
and they are reversible: by a result of Kukiela \cite{Kuk}, well founded posets with finite levels are reversible.
So, it is natural to ask what is going on with the direct products of other linear orders;
in particular, are the integer plane, $\Z ^2$, or the rational plane, $\Q ^2$, reversible lattices?
What is going on with the direct powers $\Q ^\k$ and $\Z ^\k$?

Second, Boolean algebras are reversible partial orders \cite{Kuk}
and it is natural to ask what is going on with some relevant suborders of them (e.g.\ ideals, filters).
For example, if $\o \leq \l \leq \k$ are cardinals, is the partial order $\la [\k]^{<\l }, \subset\ra$ of $<\l$-sized subsets of $\k$ reversible?
What is going on with the ideal of Borel meager subsets of a Polish space?

Third, some countable ultrahomogeneous structures are reversible (e.g.\ the rational line, Henson graphs \cite{KuMo2} and digraphs \cite{KRet})
and some are not (e.g.\ the random graph). What is going on with the random poset, or, more generally, with $\k$-universal-homogeneous posets?

In this paper we present a method for proving non-reversibility of structures based on some techniques from model theory
and detect several classes of non-reversible partial orders containing the particular posets mentioned above.
We recall that back and forth systems of partial isomorphisms between two structures of the same language, $\X$ and $\Y$,
were used in model theory for characterizations of
isomorphism of countable structures (Cantor),
elementary equivalence (Ehrenfeucht and Fra\"{\i}ss\'{e}, \cite{Ehr,Fra1,Fra2}),
$L_{\infty \o}$-equivalence (Karp \cite{Karp}).
Using these ideas we first give a simple characterization of (non-)reversibility;
roughly: an $L$-structure $\X$ is not reversible
iff there is a back and forth system (b.f.s.) $\Pi$ of partial self-condensations of $\X$
containing one which is not a partial isomorphism
and having certain closure properties.

In the rest of the paper we use the non-trivial direction of that theorem.
Roughly speaking, if we suspect that a partial order $\X$ it is not reversible,
our task is to show that a condensation $f\in \Cond (\X )\setminus \Aut (\X )$ exists;
so, the strategy is to construct a b.f.s.\ $\Pi$ consisting of approximations of the ``unknown" self-condensation $f$.

It is reasonable to expect that the approximations $\f \in \Pi$ are ``small";
for example, that $|\f |<|X|$, or that $\dom \f$ is a convex set, or that the field of $\f$, that is, the set $\dom \f \cup \ran \f$, is bounded.

In addition, we have to ensure the existence of a filter $\Phi$ in the partial order $\la \Pi ,\supset \ra$
such that $\bigcup \Phi$ is a bijection from $X$ onto $X$ and, hence, $f:=\bigcup \Phi \in \Cond (\X )\setminus \Aut (\X )$, as desired.
This can be done by requesting that some cardinal invariants of the poset $\Pi$ are large;
for example, that $\fm _\Pi >|X|$, where $\fm _\Pi$ is the Baire number of $\Pi$, or more strongly, that $\Pi$ is $|X|$-closed.
If $|X|=\o$, then the last condition is automatically satisfied by the Rasiowa-Sikorski theorem.
If $|X|$ is a singular cardinal then the last condition can be too strong;
hence in Theorem \ref{T8211} we avoid it, but put additional requirements on the b.f.s.\ $\Pi$.

While the proof of Theorem \ref{T8196} characterizing non-reversible $L$-structures is elementary,
its applications strongly depend of the properties of the (classes) of structures under consideration;
so, we have to develop different strategies in different situations.
In this paper some strategies are presented and applied in the proofs that several (classes of) {\it connected} partial orders are not reversible.
We note that the reversibility of {\it disconnected} binary structures was characterized \cite{KuMo5}
and several specific results concerning that topic can be found in \cite{KuMo4} and \cite{KuMo3}.

Finally we remark that the results concerning the non-reversibility of partial orders
can be regarded as the results concerning non-reversibility of the corresponding topological spaces:
if $\P=\la P,\leq \ra$ is a partial order and $\CO$ is the topology on the set $P$ generated by the base consisting of the proper ideals, $(\cdot ,p]$, $p\in P$,
then a mapping $f :P\rightarrow P$ is a homomorphism of $\P$ iff it is continuous with respect to $\CO$.
Thus the poset $\P$ is reversible iff $\la P, \CO\ra$ is a reversible topological space.
Also the non-reversibility of $\P$ implies the non-reversibility of the corresponding comparability graph (see Remark 5 of \cite{Kuk}).
\section{Preliminaries}\label{S2}
If $\P =\la P, \leq \ra $ is a partial order, the elements $p$ and $q$ of $P$ are {\it incompatible},
we write $p\perp q$, iff there is no $r\in P$ such that $r\leq p,q$.
A set $A\subset P$ is an antichain iff $p\perp q$, for different $p,q\in A$.
$p\in P$ is a {\it minimal element} of $\P$ iff there is no $q<p$.
$p\in P$ is an {\it atom } of $\P$ iff each two elements $q,r\leq p$ are compatible.
By $\Min (\P)$ and $\At (\P )$ we denote the sets of minimal elements and atoms respectively. Clearly, $\Min (\P)\subset \At (\P )$.

A set $D\subset P$ is called {\it dense} iff for each $p\in P$ there is $q\in D$ such that $q\leq p$;
$D$ is called {\it open} iff $q\leq p\in D$ implies $q\in D$.
It is easy to check that a set $D\subset P$ is dense (resp.\ open) iff $D$ is a dense (resp.\ open) set in the topological space $\la P, \mathcal{O}\ra$,
where $\mathcal{O}$ is the topology on the set $P$ generated by the base $\CB :=\{ (\cdot ,p] :p\in P\}$.
By $\DO (\P )$ (resp.\ $\RD(\P )$) we denote the set of all dense open (resp.\ dense) subsets of $\P$.
We note that open sets are also called {\it initial segments} or {\it downwards closed sets}
and that, for a set $S\subset P$, $S\!\downarrow :=\{ p\in P: \exists s\in S \;p\leq s\}$ is the minimal open superset of $S$ (see \cite{Harz});
thus a set $S$ is open iff  $S\!\downarrow=S$.
Dually we define {\it final segments} and the sets $S\!\uparrow$.

A set $\Phi\subset P$ is a {\it filter} in $\P$ iff $\Phi \ni p\leq q$ implies $q\in \Phi$
and for each $p,q\in \Phi$  there is $r\in \Phi$ such that $r\leq p,q$.
If $\CD \subset \RD (\P)$, then a filter  $\Phi$ is called {\it $\CD$-generic} iff $\Phi \cap D \neq \emptyset$, for each $D\in \CD$.
It is easy to check that for each family $\CD \in [\DO (\P)]^\k$ there is a $\CD$-generic filter iff the same holds for each family $\CD \in [\RD (\P)]^\k$.

A partial order $\P =\la P, \leq \ra $ is called {\it $\k$-closed} (where $\k\geq \o$ is a cardinal)
iff whenever $\g <\k$ is an ordinal and $\la p_\a : \a <\g \ra$ is a sequence in $P$ such that $p_\b \leq p_\a$, for $\a <\b<\g$,
there is $p\in P$ such that $p\leq p_\a$, for all $\a <\g$. So if $\P$ contains a chain without a lower bound
let $\ft _\P := \min \{ \k: \P \mbox{ is not } \k\mbox{-closed}\}$; otherwise, $\ft _\P:=\infty$.
\begin{fac}\label{T8192}
If $\,\P $ is a $\k$-closed partial order, then

(a) If $\CD \subset \RD (\P)$ and $|\CD|\leq \k$, then there is a  $\CD$-generic filter $\Phi$ in $\P$;

(b) If $\k$ is a singular cardinal, then $\,\P $ is $\k ^+$-closed.
\end{fac}
\dok
(a) Let $\CD =\{ D_\a :\a <\k \}$ be an enumeration.
By recursion we construct a sequence $\la p_\a :\a <\k \ra$ in $P$ such that for all $\a ,\b \in \k$ we have
(i) $p_\a \in D_\a$, and
(ii) $\a <\b \Rightarrow p_\b \leq p_\a $.
First we take $p_0\in D_0$. Suppose that $0<\a <\k$ and that $\la p_\b :\b <\a \ra$ is a sequence satisfying (i) and (ii).
Then, since $\P $ is $\k$-closed, there is $p\in P$ such that $p\leq p_\b$, for all $\b <\a$, and, since $D_\a$ is dense in $\P$,
there is $p_\a \in D_\a$, such that $p_\a \leq p$, which implies that $p_\a \leq p_\b$, for all $\b <\a$.
Thus the sequence $\la p_\b :\b \leq\a \ra$ satisfies (i) and (ii) and the recursion works.
It is evident that $\Phi :=\{ p\in P : \exists \a <\k \; p_\a \leq p\}$ is a filter  in $\P$.

(b) See \cite{Kun}, p.\ 240.
\hfill $\Box$
\section{Back-and-forth systems and reversibility of relational structures}\label{S3}
Let $L=\la R_i:i\in I\ra$ be a relational language, where $\ar (R_i)=n_i\in \N$, for $i\in I$.
In this section we give a characterization of (non-)reversible $L$-structures, which will be used in the sequel.
If $\X $ is an $L$-structure,
a function $\f$ will be called a {\it partial condensation of $\X$}, we will write $\f \in \Pc (\X )$, iff
$\f$ is a bijection which maps $\dom (\f )\subset X$ onto $\ran (\f )\subset X$ and
\begin{equation}\label{EQ0092}
\forall i\in I \;\; \forall \bar x \in (\dom (\f ))^{n_i} \;\;(\bar x \in R_i^\X \Rightarrow \f\bar x \in R_i^\X ).
\end{equation}
A set $\Pi \subset \Pc (\X )$ will be called a {\it back and forth system of condensations} (shortly: b.f.s.)
iff $ \Pi\neq \emptyset$ and 

(bf1) $\forall \f\in \Pi \;\; \forall a\in X \;\;  \exists \p \in \Pi\;\;(\f \subset \p \land a\in \dom \p )$,

(bf2) $\forall \f\in \Pi \;\; \forall b\in X \;\;\,\exists \p \in \Pi\;\;(\f \subset \p \land b\in \ran \p )$.

\noindent
We will say that a b.f.s.\ $\Pi $ is a {\it $\k$-closed b.f.s.}
iff the partial order $\la \Pi ,\supset\ra$ is $\k$-closed.
A partial condensation $\f ^* \in \Pc (\X )$ will be called {\it bad},
we will write $\f ^*\in \Bad (\X)$,
iff  $\f ^*$ is not a partial isomorphism of $\X$;
that is, there are $i\in I$ and $\bar x =\la x_0 ,\dots ,x_{n_i-1}\ra\in \dom(\f )^{n_i}$
such that $\bar x \not\in R_i^\X$ and $\f ^*\bar x \in R_i^\X$.
Then we define $\Pc _{\f ^*}(\X ):=\{ \f \in \Pc (\X ) : \f ^* \subset \f \}$.
If $\k$ is a cardinal let  $\Pc ^{<\k }(\X ):=\{ \f \in \Pc (\X ) : |\f |<\k \}$ and $\Pc _{\f ^*}^{<\k }(\X ):=\{ \f \in \Pc^{<\k } (\X ) : \f ^* \subset \f \}$.

\begin{te}\label{T8196}
For an $L$-structure $\X $ of size $\k\geq \o$ the following is equivalent:

(a) $\X$ is not a reversible structure,

(b) There are $\f ^* \in \Bad (\X)$ and a $\k$-closed b.f.s.\ $\Pi\subset \Pc _{\f ^*} (\X )$,

(c) There exists a b.f.s.\ $\Pi\subset \Pc (\X )$ containing a bad condensation, if $\k =\o$.
\end{te}
\dok
(a) $\Rightarrow$ (b).
If $f\in \Cond (\X )\setminus \Aut (\X )$,
then $f \in \Bad (\X )$,
$\Pi := \{ f \}$ is a b.f.s.\
and each chain in $\Pi$ has a lower bound;
thus $\Pi$ is $\k$-closed.

(b) $\Rightarrow$ (a).
Let $\f ^* \in \Bad (\X)$, let $\Pi\subset \Pc _{\f ^*} (\X )$ be a $\k$-closed b.f.s.\ and $\P :=\la \Pi ,\supset\ra$.
Let $A:=X\setminus \dom (\f ^*)$ and $B:=X\setminus \ran (\f ^*)$.
For $a \in A$, let $D_a :=\{ \p \in \Pi : a \in \dom \p \}$.
If $\f  \in \Pi $, then   by (bf1) there is $\p \in \Pi$ such that $a\in \dom \p$ and $\f  \subset \p $;
thus, $\p \in D_a$ and $\p \supset  \f $.
So, the sets $D_a$, $a\in A$, are dense in $\P$
and, similarly, the sets $\Delta _b :=\{ \p \in \Pi  : b \in \ran \p \}$, $b\in B$, are dense in $\P$.
Since $|A|+|B|\leq \k$, by Fact \ref{T8192}(a) there is a filter $\Phi $ in $\P$ intersecting all $D_a$'s
and $\Delta _b$'s. Clearly we have $\f ^* \subset F:= \bigcup \Phi\subset X \times X $.

If $\la x,y'\ra , \la x,y''\ra \in F$,
there are $\f ' ,\f '' \in \Phi$ such that $\la x,y'\ra \in \f'$ and $\la x,y''\ra \in \f ''$
and, since $\Phi$ is a filter, there is $\f \in \Phi$ such that $\f \supset \f ' , \f ''$.
Thus $\la x,y'\ra , \la x,y''\ra \in \f$
and, since $\f$ is a function, we have $y'=y''$.
So $F$ is a function and in the same way we show that it is an injection.
If $a\in X\setminus \dom (\f^*)$,
then there is $\f\in \Phi \cap D_a$
and, hence, $a\in \dom (\f ) \subset \dom (F)$.
So $X\setminus \dom (\f^*)\subset \dom (F)$,
which, together with $\dom (\f ^*)\subset \dom (F)$, implies that $\dom (F)=X$.
Similarly we have $\ran (F)=X$ and, thus, $F :X \rightarrow X$ is a bijection.

Let $i\in I$, $\bar x =\la x_0, \dots ,x_{n_i-1}\ra\in (\dom (F) )^{n_i}$ and $\bar x \in R_i^\X$.
Then, since $\dom (F)=\bigcup _{\f \in \Phi}\dom (\f )$,
for each $j<n_i$ there is $\f _j\in \Phi$ such that $x_j \in \dom (\f _j)$.
Since $\Phi$ is a filter, there is $\f \in \Phi$ such that for each $j<n_i$ we have $\f \supset \f _j$
and, hence, $\dom (\f ) \supset \dom (\f  _j)$.
Thus $\bar x\in (\dom (\f ) )^{n_i}$
and, since $\f \in \Pc (\X )$ and $\bar x \in R_i^\X$, we have $F\bar x =\f \bar x\in R_i^\X$.
So $F\in \Cond (\X )$ and, since $\f ^*\subset F$, we have $F\not\in \Aut (\X )$.
Thus, the structure $\X $ is not reversible.

The equivalence (b) $\Leftrightarrow$ (c) is true because each poset is $\o$-closed.
\hfill $\Box$
\begin{rem}\rm
Condition (b) of Theorem \ref{T8196} can be written as
\begin{equation}\label{EQ8183}
\exists \f ^* \in \Bad (\X) \;\;\exists \mbox{b.f.s.} \Pi\subset \Pc _{\f ^*} (\X )\;\; \ft _\Pi >|X|.
\end{equation}
In fact,  $\ft _\Pi$ can be equivalently replaced by some another cardinal invariants of the poset $\la \Pi ,\supset\ra$.
Namely, if $\P$ is a partial order, then there is a family $\CD \subset \DO (\P )$ such that $\bigcap \CD $ is not dense (resp.\ $\bigcap \CD =\emptyset$, there is no $\CD$-generic filter in $\P$)
iff $\Min (\P)$ is not dense (resp.\ $\Min (\P)=\emptyset$, $\At (\P)=\emptyset$), see \cite{Kun}. Then
\begin{eqnarray*}
\textstyle
\add (\CN  _\P) & := & \textstyle\min \{ |\CD|: \CD \subset \DO (\P ) \mbox{ and }\bigcap \CD \mbox{ is not dense}\},\\
\textstyle
\cov (\CN  _\P) & := & \textstyle\min \{ |\CD|: \CD \subset \DO (\P ) \mbox{ and }\bigcap \CD =\emptyset\}, \\
\textstyle
\fm _\P         & := &\textstyle \min \{ |\CD|: \CD \subset \DO (\P ) \mbox{ and there is no $\CD$-generic filter in $\P$}\};
\end{eqnarray*}
otherwise, these invariants are defined to be $\infty$.\footnote{
If $D\in \DO (\P )$, then $P\setminus D$ is a closed nowhere dense set. Let $\CN _\P$ denote the ideal of nowhere dense subsets of $\P$.
The equalities $\add (\CN  _\P)=\min \{ |\CN |: \CN \subset \CN_\P  \mbox{ and }\Int (\bigcup \CN )\neq \emptyset \}$
and $\cov (\CN  _\P)= \min \{ |\CN|: \CN \subset \CN _\P \mbox{ and }\bigcup \CN =P \}$
explain the notation $\add (\CN  _\P)$ and $\cov (\CN  _\P)$ (when they are $<\infty$).}
Now, since the inequality $\fm _\Pi>|X|=:\k$ provides the filter intersecting $\k$-many dense sets $D_a$ and $\Delta _b$,
in (\ref{EQ8183}) we can replace  $\ft _\Pi>|X|$ by $\fm _\Pi>|X|$.
Finally, it is known that $\add (\CN  _\P) \leq  \cov (\CN  _\P)\leq \fm _\P$ (see \cite{Bart})
so, in (\ref{EQ8183}) we can replace  $\ft _\Pi>|X|$ by $\add (\CN  _\P)>|X|$ or by $\cov (\CN  _\P)>|X|$,
because these conditions imply $\fm _\Pi>|X|$.

In addition, some forcing axioms can help in this context.
For example, Martin's axiom (MA) says that for each ccc partial order $\P$ we have $\fm _\P \geq \fc$ (see \cite{Kun}).
Thus in models of MA we have: an $L$-structure $\X $ of size $<\fc$ is not reversible  iff there is a ccc b.f.s.\ $\Pi\subset \Pc (\X )$ containing a bad condensation.
Similarly, the Proper forcing axiom (PFA) says that for each proper partial order $\P$ we have $\fm _\P \geq \fc (=\o _2)$;
so, in models of PFA we have: an $L$-structure $\X $ of size $\o_1$ is not reversible  iff there is a proper b.f.s.\ $\Pi\subset \Pc (\X )$ containing a bad condensation.
\end{rem}

We intend to use Theorem \ref{T8196}
in proofs that some structures are not reversible,
by construction of the convenient $\k$-closed b.f.s.\ $\Pi$
of ``small" approximations of the desired condensation $f\in \Cond (\X )\setminus \Aut (\X )$.
Roughly speaking, if the cardinal $\k$ is singular, then by Fact \ref{T8192}(b) the b.f.s.\ $\Pi$ will be $\k^+$-closed and it can happen that $f\in \Pi$ although $f$ is not ``small",
which is impossible (see Example \ref{EX8100}).
Theorem \ref{T8211} given in the sequel provides sufficient conditions for non-reversibility of structures of singular cardinality and the assumption that $\Pi$ is $\k$-closed is omitted.
Let $X$ be a set of size $\k\geq \o$. We will say that a mapping $\cl : P(X)\rightarrow P(X)$ is a {\it closure operator on $P(X)$} iff

(cl1) $S\subset \cl (S)$, for each $S\subset X$,

(cl2) $\cl (\bigcup _{i\in I}S_i)= \bigcup _{i\in I}\cl (S_i)$, whenever $|I|<\k$ and $|S_i|<\k$, for all $i\in I$.

\noindent
Also we will regard the following strengthening of conditions (bf1) and (bf2)

(bf1$_\l$) $\forall \f\in \Pi \;\; \forall a\in X \;\;\exists \p \in \Pi\;\;(\f \subset \p \land a\in \dom \p \land |\p \setminus \f|<\l )$,

(bf2$_\l$) $\forall \f\in \Pi \;\; \forall b\in X \;\;\,\exists \p \in \Pi\;\;(\f \subset \p \land b\in \ran \p \land\,|\p \setminus \f|<\l  )$.
\begin{te}\label{T8211}
Let $\X $ be an $L$-structure of size $\k$ and $\l$ a regular cardinal, where $\o\leq \l\leq \k $.
Let $\cl $ be a closure operator on $P(X)$ and $\f ^* \in \Pc ^{<\l } (\X )$ a bad condensation such that $\cl (\dom (\f ^*))=\dom (\f ^*)$.
If the set
\begin{equation}\label{EQ8171}
\Pi := \Big\{ \f \in \Pc _{\f ^*}^{<\k }(\X ) : \cl (\dom (\f ))=\dom (\f )\Big\}
\end{equation}
satisfies conditions (bf1$_\l$) and (bf2$_\l$), then $\X$ is not reversible.
\end{te}
\dok
Since $|\f ^*|<\l <\k$ and $\cl (\dom (\f ^*))=\dom (\f ^*)$ we have $\f ^*\in \Pi$.

If $\k$ is a regular cardinal, $\g<\k$ and $\f _\xi \in \Pi$, for $\xi <\g$, where $\f_0 \subset \f_1 \subset \dots$,
then $\f ^*\subset \f :=\bigcup _{\xi <\g }\f _\xi\in \Pc (\X )$.
By (cl2) we have $\cl (\dom (\f ) ) =\cl (\bigcup _{\xi <\g }\dom (\f _\xi ))=\bigcup _{\xi <\g }\cl (\dom (\f _\xi ))=\bigcup _{\xi <\g }\dom (\f _\xi ) =\dom (\f )$.
By the regularity of $\k$ we have $|\f | =\sum _{\xi <\g }|\f _\xi|<\k$.
Thus $\f\in \Pi$ and $\Pi$ is $\k$-closed.
By (bf1$_\l$) and (bf2$_\l$) $\Pi$ is a b.f.s.\ and by Theorem \ref{T8196} the structure $\X $ is not reversible.

If $\k$ is a singular cardinal, then by our assumptions we have $\l <\k$.
Also, $|\dom \f^*|=|\ran \f ^*|=|\f ^*|<\l<\k$ and, hence,  $|X \setminus \dom \f ^*|=|X \setminus \ran \f ^*|=\k$.
Let $E$ and $O$ be the sets of even and odd ordinals $<\k $ respectively
and let $X \setminus \dom \f ^* =\{ a _\xi :\xi \in E\}$ and $X \setminus \ran \f ^* =\{ b _\xi :\xi \in O\}$ be enumerations.
By recursion we define a sequence $\la \p _\xi :\xi <\k \ra$ such that for each $\xi <\k $ we have:

(i) $\p _\xi \in \Pi$,

(ii) $\zeta < \xi \Rightarrow \p _\zeta \subset \p _\xi$,

(iii) $a_\xi \in \dom \p _\xi $, if $\xi \in E$,

(iv) $b_\xi \in \ran \p _\xi $, if $\xi \in O$,

(v) $|\p _\xi |\leq |\xi| +\l $.

\noindent
We show that the recursion works.
By (bf1$_\l$), for $\f ^*$ and $a_0$
there is $\p _0\in \Pi$ such that $\f ^* \subset \p _0$, $a_0 \in \dom \p _0$ and $|\p _0 \setminus \f^* |<\l$;
so, (i) and (iii) are true and (ii) and (iv) are true trivially.
Since $|\p _0|=|\f ^*|+|\p _0 \setminus \f^*|<\l$ we have (v).

Let $0<\xi <\k $ and suppose that $\la \p _\zeta : \zeta <\xi \ra$ is a sequence satisfying (i)--(v).

{\it Case 1}: $\xi \in O$.
Then $\p _{\xi -1}\in \Pi$
and by (bf2$_\l$) we can choose $\p _\xi\in \Pi$
such that $\p _{\xi -1}\subset \p _\xi$, $b_\xi \in \ran \p _\xi$ and $|\p _\xi \setminus \p _{\xi -1} |<\l$;
so, (i)--(iv) are true.
If $\l \geq |\xi -1|$,
then $|\p _{\xi -1} |\leq |\xi -1|+\l =\l $
and, since $|\p _\xi \setminus \p _{\xi -1} |<\l$ and $\p _{\xi -1}\subset \p _\xi$,
we have $|\p _\xi |= |\p _{\xi -1} \cup (\p _\xi \setminus \p _{\xi -1})|= |\p _{\xi -1}| + |\p _\xi  \setminus \p _{\xi -1}|\leq \l =|\xi|+\l$
and (v) is true.
If $\l < |\xi -1|$,
then $|\p _{\xi -1} |\leq |\xi -1|+\l =|\xi -1|=|\xi| $
and $|\p _\xi \setminus \p _{\xi -1} |<\l< |\xi -1|$,
and, hence, $|\p _\xi |= |\p _{\xi -1}| + |\p _\xi  \setminus \p _{\xi -1}|\leq|\xi | =|\xi|+\l$.
Thus (v) is true again.

{\it Case 2}: $\xi \in E$ and  $\xi \not\in \Lim$.
Then $\p _{\xi -1}\in \Pi$
and by (bf1$_\l$) there is $\p _\xi\in \Pi$
such that $\p _{\xi -1}\subset \p _\xi$,
$a_\xi \in \dom \p _\xi$
and $|\p _\xi \setminus \p _{\xi -1} |<\l$; so, (i)--(iv) are true.
The proof of (v) is exactly as in Case 1.

{\it Case 3}: $\xi \in \Lim$.
Then, by (i) and (ii), $\la \p _\zeta : \zeta <\xi \ra$ is a chain of partial condensations
and, hence $\f ^*\subset \eta :=\bigcup _{\zeta <\xi }\p _\zeta\in \Pc (\X )$.
By (i) and (cl2) we have $\cl (\dom (\eta ) ) =\cl (\bigcup _{\zeta <\xi }\dom (\p _\zeta ))=\bigcup _{\zeta <\xi }\cl (\dom (\p _\zeta ))=\bigcup _{\zeta <\xi }\dom (\p _\zeta ) =\dom (\eta)$.

If $|\xi|<\l$,
then by (v) for each $\zeta <\xi$ we have $|\p _\zeta | \leq |\zeta | +\l=\l$
and, by the regularity of $\l$, $|\eta | \leq \sum _{\zeta <\xi }|\p _\zeta |\leq \l <\k $,
which gives $\eta \in \Pi$.
So, by (bf1$_\l$) there is $\p _\xi\in \Pi$ such that $\eta \subset \p _\xi$,
$a_\xi \in \dom \p _\xi$ and $|\p _\xi \setminus \eta |<\l$; thus, (i)--(iv) are true.
Since $\eta\subset \p _\xi$ and $|\p _\xi \setminus \eta |<\l$,
it follows that $|\p _\xi |= |\eta | + |\p _\xi \setminus  \eta |\leq \l= |\xi|+\l$,
and (v) is true too.

If $\l \leq |\xi|$,
then for each $\zeta <\xi$ we have $|\p _\zeta |\leq |\zeta| +\l \leq |\xi | +\l =|\xi |$
and, hence, $|\eta | \leq \sum _{\zeta <\xi }|\p _\zeta |\leq |\xi | <\k $,
which gives $\eta \in \Pi$.
By (bf1$_\l$) there is $\p _\xi\in \Pi$ such that $\eta \subset \p _\xi$,
$a_\xi \in \dom \p _\xi$ and $|\p _\xi \setminus \eta |<\l$; so, (i)--(iv) are true.
Since $\eta\subset \p _\xi$ and $|\p _\xi \setminus \eta |<\l$,
it follows that $|\p _\xi |= |\eta | + |\p _\xi \setminus  \eta |\leq  |\xi|+\l$,
and (v) is true.
Thus, the recursion works.

Now, by (i), (ii)  and the definition of $\Pi$
we have $\f ^* \subset f:= \bigcup _{\xi <\k} \p _\xi \in \Pc (\X )$,
by (iii) and (iv) it follows that $\dom f =\ran f =X$,
which means that $f\in \Cond (\X )$.
Since $\f ^* \subset f$ we have $f\not\in \Aut (\X )$ and $\X$ is not reversible.
\hfill $\Box$
\section{Small approximations. Homogeneous-universal posets}\label{S4}
Since the class of posets is a J\'{o}nsson class \cite{Jon}, for each regular beth number $\k$ there is a $\k$-homogeneous-universal poset $\P$.
As an example of application of Theorem \ref{T8196} we show that such posets are not reversible; in particular, taking $\k =\o$ we conclude that  the random poset
(i.e., the unique countable homogeneous universal poset, see \cite{Sch}) is non-reversible as well.
If $\mathbb P$ is a poset and $p,q \in P$ are incomparable elements (i.e.\ $p\not\leq q$ and $q\not\leq p$), we will write $p\parallel q$. For $A,B\subset P$,
$A<B$ denotes that $a<b$, for all $a\in A$ and $b\in B$; notation $A\parallel B$ is defined similarly. We regard the posets of size $\k \geq \o $ satisfying the following conditions

{\rm (u1)} $\forall L,G \in [P ]^{<\k }\setminus \{ \emptyset \} \;\; (L<G \Rightarrow \exists p \in P \;\;L<p<G)$,

{\rm (u2)} $\forall K \in [P ]^{<\k }\setminus \{ \emptyset \} \;\; \exists p,q,r \in P \;\;(p<K \land q> K  \land r\parallel  K )$.
\begin{te}\label{T8105}
If $\P= \langle P, < \rangle$ is a strict partial order of regular cardinality $\k \geq \o$
and satisfies (u1) and (u2),
then $\P$ is not reversible.
\end{te}
\dok
First we show that for each $L,G,K \in [P ]^{<\k }\setminus \{ \emptyset \}$, where $L<G$, we have
\begin{equation}\label{EQ8181}
|\{ p \in P : L<p<G \}|=\k,
\end{equation}
\begin{equation}\label{EQ8182}
|\{ p \in P : p<K \}|=|\{ p \in P : p>K \}|=|\{ p \in P : p\parallel K \}|=\k.
\end{equation}
If $L,G \in [P ]^{<\k }\setminus \{ \emptyset \}$, $L<G$ and $S:=\{ p \in P : L<p<G \}$,
then by (u1) we have $S\neq \emptyset$
and, clearly, $L<S$.
Assuming that $|S|<\k$,
by (u1) there would be $q\in P$ such that $L<q<S$
and, hence, $q\not\in S$ and $L<q <G$.
But then $q\in S$ and we have a contradiction. So, $|S|=\k$ and (\ref{EQ8181}) is true.

If $K \in [P ]^{<\k }\setminus \{ \emptyset \}$,
then by (u2) we have  $T:=\{ p \in P : p<K \}\neq \emptyset$.
Assuming that $|T|<\k$, by (u2) there would be $p\in P$ such that $p<T$,
and, hence, $p\not\in T$ and $p<K$.
But then $p\in T$ and we have a contradiction.
So, $|T|=\k$. Similarly, $|\{ p \in P : p>K \}|=|\{ p \in P : p\parallel K \}|=\k$ and (\ref{EQ8182}) is true.

By (u2) there are $a_0,a_1,b_0,b_1\in P$, where $a_0 \parallel  a_1$ and $b_0 <b_1$
and, hence, $\f ^*:=\{ \la a_0,b_0\ra ,\la a_1,b_1\ra \}\in \Bad (\P )$.
Since $\k$ is a regular cardinal
the poset $\la \Pi \supset \ra$, where $\Pi :=\Pc _{\f ^*}^{<\k}(\P )$, is $\k$-closed
and, by Theorem \ref{T8196}, it remains to be shown that $\Pi$ is a b.f.s.

(bf1) If $\f \in \Pi$ and $a\in P \setminus  \dom \f$,
then the sets $L_a:=\{p\in \dom \f : p<a\}$ and  $G_a:=\{q\in \dom \f : q>a\}$ are of size $<\k$
and we have the following cases.

{\it Case 1}: $L_a \neq \emptyset$ and $G_a \neq \emptyset$.
If $l\in \f [L_a]$ and $g\in \f [G_a]$,
then there are $p\in L_a$ and $q\in G_a$ such that $l=\f (p)$ and $g=\f(q)$
and, since $p<a<q$ and $\f$ is a homomorphism, $\f (p)<\f (q)$, that is $l<g$.
Thus $\f [L_a] < \f [G_a] $
and, since $|\ran \f|<\k$,
by (\ref{EQ8181}) there is $b\in P \setminus \ran \f$ such that $\f [L_a] <b< \f [G_a] $.
Now $\p := \f \cup \{ \la a,b \ra\}$ is an injection, $a\in \dom \p$ and $\f \subset \p$.
If $p\in \dom \f$ and $p<a$,
then $p\in L_a$
and, hence, $\f (p)\in \f [L_a]$ and $\p(p)=\f(p)<b =\p (a)$.
Similarly, $a<q\in \dom \f$ implies $\p (a)< \p(q)$.
So $\p$ is a homomorphism and $\p \in \Pi$.

{\it Case 2}: $L_a = \emptyset$ and $G_a \neq \emptyset$.
Then $\f  [G_a]\in [P ]^{<\k }\setminus \{ \emptyset \}$
and, since $|\ran \f|<\k$, by (\ref{EQ8182}) there is $b\in P \setminus \ran \f$ such that $b< \f [G_a] $.
Now $\p := \f \cup \{ \la a,b \ra\}$ is an injection, $a\in \dom \p $, $\f \subset \p$
and we show that $\p$ is a homomorphism.
If $a<q\in \dom \f$, then $q\in G_a$
and, hence, $\f (q)\in \f [G_a]$ and $\p (a)=b <\f (q)=\p (q)$.

{\it Case 3}: $L_a \neq \emptyset$ and $G_a = \emptyset$. This case is the dual of Case 2.

{\it Case 4}: $L_a = \emptyset$ and $G_a = \emptyset$.
Then $a\parallel  \dom \f$
and choosing $b\in P\setminus \ran \f$
we have $\p := \f \cup \{ \la a,b \ra\}\in \Pi$.

(bf2) Let $\f \in \Pi$ and $b\in P \setminus  \ran \f$.
Since $|\dom \f |<\k$, by (\ref{EQ8182}) there is $a\parallel  \dom \f$.
Thus $\p := \f \cup \{ \la a,b \ra\} \in \Pi $ and (bf2) is true.
\hfill $\Box$
\begin{prop}\label{T8219}
(a) If $\,\P$ is a $\k$-homogeneous-universal poset, it is not reversible.

(b) The random poset is not reversible.
\end{prop}
\dok
(a) Let $\P$ be a $\k$-homogeneous-universal poset (of regular size $\k$).
Thus, each isomorphism between $<\k$-sized substructures of $\P$ extends to an automorphism of $\P$
and each poset of size $\leq \k$ embeds in $\P$.
By Theorem \ref{T8105} it is sufficient to show that $\P$ satisfies (u1) and (u2).

Let $L,G \in [P ]^{<\k }\setminus \{ \emptyset \}$, where $L<G$,
and let $\Y$ be a poset with domain $Y=Y_L \cup \{ a\} \cup Y_G$,
where $Y _L < \{ a\}< Y _G$, $\Y _L\cong \BL$ and  $\Y _G\cong \BG$.
Then $|Y|<\k$ and, by the universality of $\P$ there is an embedding $e:\Y \hookrightarrow \P$.
Thus $L\cong e[Y_L]<e(a)<e[Y_G]\cong G$
and there are isomorphisms $f_L:e[Y_L]\rightarrow L$ and $f_G:e[Y_G]\rightarrow G$.
Clearly, $f:=f_L \cup f_G :e[Y_L]\cup e[Y_G]\rightarrow L\cup G$ is an isomorphism between $<\k$-sized substructures of $\P$
and, by the homogeneity of $\P$, there is an automorphism $F\in \Aut (\P )$ such that $f\subset F$,
which implies that $L <p:=F(e(a))<G$.
So (u1) is true and (u2) has a similar proof.

(b) It is known (see \cite{Cam}) that a countable poset $\mathbb P $ is isomorphic to the random poset iff
for each triple $\langle L,G,U \rangle$ of pairwise disjoint elements of $[P]^{<\o}$ satisfying
(C1) $L < G $,
(C2) $\forall u\in U \; \forall l\in L \;\; \neg u<l $,
and
(C3) $\forall u\in U \; \forall g\in G \;\; \neg g<u $,
there exists $p\in P$ such that
$L < p <G $ and $ p \parallel u $, for all $u\in U$.
So, the random poset satisfies (u1) and (u2) and, by Theorem \ref{T8105}, it is not reversible.
\hfill $\Box$
\section{Small approximations with open domain.\\ Well founded posets}\label{S5}
A poset $\P =\la P,\leq\ra$ is {\it well founded} iff $\o ^* \not\hookrightarrow \P$, iff each non-empty set $A\subset P$ has a minimal element.
Then (see \cite{Harz}, p.\ 231) there is an ordinal $\xi$ (the {\it height of $\,\P$}, $\he \P$) such that $P=\bigcup _{\a < \xi }L_\a$,
where the {\it levels of $\,\P$} are defined by $L_0=\Min (P)$ and
$L_\a=\Min (P\setminus \bigcup _{\b <\a}L_\b)$.
$\{ L_\a :\a <\xi \}$ is a partition of $P$ into antichains
and the {\it height of $p\in P$} is defined by $\he (p)=\a$ iff $t\in L_a$.
\begin{fac}\label{T8223}
If $\P $ is a well founded poset and $p,q \in P$, then (see \cite{Harz}, p.\ 232, 248)

(a) $p<q \Rightarrow \he p < \he q$;

(b) $\a <\he (p) \Rightarrow (\cdot ,p] \cap L_\a \neq \emptyset$;

(c) There is a well order $\vartriangleleft$ on the set $P$ which extends $<$ (i.e., $<\;\subset \,\vartriangleleft$).
\end{fac}
We note that the reversible well founded posets representable as disjoint unions of well orders are characterized in \cite{KuMo4}.
While such posets are disconnected, Theorem \ref{T8100} given in the sequel is related to directed and, hence, connected ones.
In its proof we will use the following consequence of Theorem \ref{T8211}.
\begin{te}\label{T8118}
Let $\P $ be a poset of size $\k$,
$\l$ a regular cardinal, where $\o \leq \l\leq \k $,
and $\f ^* \in \Pc^{<\l } (\P )$ a bad condensation such that $\dom \f ^*$ is an open set.
If the set
$$
\Pi :=\Big\{ \f \in \Pc _{\f ^*}^{<\k }(\P ) :  (\dom \f ) \!\downarrow \;= \dom \f \Big\}
$$
satisfies conditions (bf1$_\l$) and (bf2$_\l$), then $\P$ is not reversible.
\end{te}
\dok
For $S\subset P$ let $\cl (S)= S\!\downarrow =\{ p\in P : \exists s\in S \; p\leq s\}$.
Clearly, $\cl$ is a closure operator on $P(P)$
and we apply Theorem \ref{T8211}.
\kdok
\noindent
If $\l$ is a cardinal, a poset $\P$ is called {\it $<\l$-directed} iff each set $S\subset P$ of size $<\l$ has an upper bound.
The least element of $\P$, if it exists, will be called the {\it root} of $\P$.
\begin{te}\label{T8100}
If $\P$ is a rooted well founded $<\l$-directed poset, $|P|=|L_1|=\k \geq \l=\he (\P)\geq \o $, where $\l $ is a regular cardinal,
and $|(\cdot ,p]|<\l$, for all $p\in P$, then $\P$ is not reversible.
\end{te}
\dok
First we show that
\begin{equation}\label{EQ8172}
\forall p\in P \;\; \forall S\in [P]^{<\k}\;\;|S\!\downarrow \cup (\cdot ,p]|<\k).
\end{equation}
Let $p\in P$ and $S=\{ p_\xi :\xi <\mu\}\subset P$, where $|S|=\mu <\k$.
Then $|S\!\downarrow|=|\bigcup _{\xi <\mu}(\cdot ,p_\xi] |\leq \sum _{\xi <\mu}| (\cdot ,p_\xi]|$ and $| (\cdot ,p_\xi]|<\l$, for all $\xi <\mu$.
If $\l=\k$, then by the regularity of $\k$ we have $|S\!\downarrow|<\k$.
If $\l<\k$, then $|S\!\downarrow|\leq \mu \l <\k$ again.
Thus, since $|S\!\downarrow|<\k$ and $| (\cdot ,p]|<\l\leq \k$, we have $|S\!\downarrow \cup (\cdot ,p]|<\k$. Thus (\ref{EQ8172}) is true.

Second, we show that
\begin{equation}\label{EQ8173}
\forall p\in P \; \forall S\in [P]^{<\k}\;\forall \zeta <\l \;\exists  \la p_\xi :\xi <\zeta \ra \in ((p,\cdot )\setminus S)^\zeta \;(\la p_\xi \ra \mbox{ strictly increasing}).
\end{equation}
Let $p\in P$, $S\in [P]^{<\k}$ and $\zeta \in\l$.
By recursion we construct sequences $\la p_\xi :\xi <\zeta\ra$ and $\la a_\xi :\xi <\zeta\ra$ such that for each $\xi <\zeta$

(i) $p<p_0$,

(ii) $a_\xi \in L_1 \setminus (S\!\downarrow \cup \bigcup _{\xi ' <\xi}(\cdot , p_{\xi '}])$

(iii) $p_\xi \geq a_\xi$,

(iv) $\forall \xi ' <\xi \;\; p_{\xi '} <p_\xi$,

(v) $p_\xi \not\in S$.

\noindent
We show that the recursion works.
First let $\xi =0$. Since $|L_1|=\k$, by (\ref{EQ8172}) there is $a_0 \in L_1 \setminus (S\!\downarrow \cup (\cdot , p])$ and (ii) is true.
Since $\P$ is $<\l$-directed there is $p_0\geq a_0,p$; so (iii) is true and (iv) is true trivially.
Assuming that $p\not<p_0$ we would have $p=p_0 \geq a_0$ and, hence $a_0\in (\cdot , p]$, which is false; thus (i) is true.
Assuming that $p_0\in S$ we would have $a_0 \in S\!\downarrow$, which is false; thus we have (v).

Let $0<\xi <\zeta$ and let   $\la p_{\xi '} :\xi ' <\xi\ra$ and $\la a_{\xi '} :\xi ' <\xi\ra$ be sequences satisfying (i)--(v).

If $\xi =\xi ^* +1$, then by (\ref{EQ8172}) there is $a_\xi \in L_1 \setminus (S\!\downarrow \cup (\cdot , p_{\xi ^*}])$.
By (iv) we have $\bigcup _{\xi '<\xi}(\cdot , p_{\xi '}]=(\cdot , p_{\xi ^*}] $ and (ii) is true.
Since $\P$ is $<\l$-directed there is $p_\xi\geq a_\xi,p_{\xi ^*}$; so (iii) is true.
Assuming that $p_{\xi ^*}\not< p_0$ we would have $p_{\xi ^*}=p_\xi \geq a_\xi$ and, hence $a_\xi\in (\cdot , p_{\xi ^*}]$, which is false; thus (iv) is true.
$p_\xi\in S$ would imply that $a_\xi \in S\!\downarrow$, which is false; thus we have (v).

If $\xi $ is a limit ordinal, then, since $|\{ p_{\xi '} :\xi ' <\xi\}|=|\xi|<\l$ and $\P$ is $<\l$-directed,
there is $q\in P$ such that $p_{\xi '}<q$, for all $\xi ' <\xi$.
By (\ref{EQ8172}) there is $a_\xi \in L_1 \setminus (S\!\downarrow \cup (\cdot , q])$.
By (iv) we have $\bigcup _{\xi '<\xi}(\cdot , p_{\xi '}] \subset(\cdot , q] $ and (ii) is true.
Since $\P$ is $<\l$-directed there is $p_\xi\geq a_\xi,q$; so (iii) is true.
For $\xi ' <\xi$ we have $ p_{\xi '} <q \leq p_\xi$; thus (iv) is true.
Assuming that $p_\xi\in S$ we would have $a_\xi \in S\!\downarrow$, which is false; thus we have (v).

Thus the recursion works. By (i), (iv) and (v) $\la p_\xi :\xi <\zeta \ra$ is a strictly increasing sequence in $(p,\cdot )\setminus S$.
So, (\ref{EQ8173}) is true.

Let $L_0=\{ r\}$ and $a_0,a_1 \in L_1$, where $a_0\neq a_1$.
By (\ref{EQ8173}) we have $(a_0 ,\cdot )\neq \emptyset$
and since $\P$ is well founded we can choose a minimal element $b_0$ of $(a_0 ,\cdot )$.
Clearly the mapping $\f ^* =\{ \la r ,r \ra,\la a_0,b_0\ra,\la a_1, a_0 \ra\} \in \Pc (\P )$ is a bad condensation
and the set $\dom \f ^* =\{ r , a_0, a_1\}$ is open. Thus $\f ^* \in \Pi$, where
\begin{equation}\label{EQ8158}\textstyle
\Pi :=\big\{ \f \in \Pc _{\f ^*}^{<\k }(\P ) :  (\dom \f ) \!\downarrow \;= \dom \f \big\}
\end{equation}
and in order to apply Theorem \ref{T8118} we show that (bf1$_\l$) and (bf2$_\l$) are true.

(bf1$_\l$) Let $\f \in \Pi$ and $a\not\in \dom \f$.
Then, since $\f ^* \subset \f$ we have $r\in \dom \f$ and, hence,  $a\neq r$.
Since the set $\dom \f $ is open, the set $D=\dom \f \cup (\cdot , a]$ is open too.
Since $|(\cdot ,a]|<\l \leq \k$ and $|\dom \f|=|\f|<\k$ we have $|D|<\k$.

Clearly $D':=(\cdot ,a]\setminus \dom \f$ is a well founded suborder of $\P$,
by Fact \ref{T8223}(c) there is a well order $\vartriangleleft$ on the set $D'$
which extends the initial order $<\upharpoonright D'$
and, hence, there is an ordinal $\zeta \cong \la D' ,\vartriangleleft\ra$.
Since $|D'|<\l$ we have $\zeta <\l$.

Since $|\dom \f \cap (\cdot ,a ]|<\l$
it follows that $|\f [\dom \f \cap (\cdot ,a ]]|<\l$
and, since $\P$ is $<\l$-directed, there is $p\in P$ such that
\begin{equation}\label{EQ8161}\textstyle
\f \big[\dom \f \cap (\cdot ,a ]\big] \leq p.
\end{equation}
Since $|\ran \f|<\k$,
by (\ref{EQ8173}) there is a strictly increasing sequence $\la p_\xi :\xi <\zeta \ra$ in $(p,\cdot )\setminus \ran \f$.
So we obtain a chain $M=\{ p_\xi : \xi <\zeta \}\cong \zeta$ in $(p,\cdot ) \setminus \ran \f $.

Let $\p :D\rightarrow P$, where $\p \upharpoonright \dom \f =\f$
and $\p \upharpoonright D' : \la D' ,\vartriangleleft \ra \rightarrow \la M,<\ra$ is an isomorphism.
Since $M \cap \ran \f=\emptyset$ the mapping $\p$ is an injection.

For a proof that $\p$ is a homomorphism
we take $x,y\in D$, where $x<y$, and prove that $\p (x)< \p (y)$.
Since the set $\dom \f$ is open it is impossible that $x\not\in \dom \f \ni y$.
If $x,y\in \dom \f $,
then $\p (x)=\f (x)< \f (y)=\p (y)$, because $\f \in \Pc (\P )$.

If $x,y\not\in \dom \f $,
then $x,y\in D'$ and, hence, $x\vartriangleleft y$,
which gives $\p (x)< \p (y)$.

If $x\in \dom \f \not\ni y$,
then $y\in (\cdot ,a]$ and, hence $x<a$,
which implies that $x\in \dom \f \cap (\cdot ,a ]$
and $\f (x) \in \f [\dom \f \cap (\cdot ,a ]]$;
so, by (\ref{EQ8161}), $\f (x)\leq p$.
In addition,  $y\in D'$ and, hence, $\p (y)\in M\subset (p, \cdot )$,
which gives $\p (x)=\f (x)\leq p < \p (y)$.
So, $\p$ is a homomorphism indeed.

Thus $\p \in \Pc (\P)$, $\dom \p =D$ is an open set, $|\p|=|D|<\k$, $\f ^*\subset \p $ and by (\ref{EQ8158}) we have $\p \in \Pi$.
In addition, $\f \subset \p $, $a\in \dom \p$ and $|\p \setminus \f|=|(\cdot ,a]\setminus \dom \f|<\l$; so (bf1$_\l$) is true.

(bf2$_\l$) Let $\f \in \Pi$ and $b\not\in \ran \f$.
Since $|\dom \f|=|\f|<\k $ there is $a\in L_1 \setminus \dom \f$ and $\p := \f \cup \{ \la a ,b \ra\}$ is an injection.
Since $r\in \dom \f$ the set $\dom \p = \dom \f \cup\{ a\}$ is open.

For a proof that $\p$ is a homomorphism, assuming that $x,y\in \{ a\} \cup \dom \f$ and $x<y$ we show that $\p (x)< \p (y)$.
Since the set $\dom \f$ is open it is impossible that $x\not\in \dom \f \ni y$
and since $\dom \p \setminus \dom \f =\{ a\}$, $x,y\not\in \dom \f $ is impossible too.
If $x,y\in \dom \f $, then, again, $\p (x)=\f (x)< \f (y)=\p (y)$.
If $x\in \dom \f \not\ni y$, then $y=a$ and, hence $x=r$.
So, since $b \not\in \ran \f \supset \ran \f ^*$ and, hence, $b\neq r$,
we have $\p (x)=\f^* (r) = r < b= \p (a)=\p (y)$.

Thus $\p \in \Pc (\P)$, $\dom \p =\dom \f \cup\{ a\}$ is an open set, $|\p|<\k$, $\f ^*\subset \p $ and by (\ref{EQ8158}) we have $\p \in \Pi$.
In addition, $\f \subset \p $, $b\in \ran \p$ and $|\p \setminus \f|=1<\l$; so (bf2$_\l$) is true.
\kdok
It is evident that if a partial order $\P$ has a root $r$, then $\P$ is reversible iff $\P \setminus \{ r\}$ is reversible.
Thus, Theorem \ref{T8100} has an evident ``rootless" modification. In addition, these two variants of Theorem \ref{T8100} have their obvious dual versions.
\begin{prop}\label{T8213}
For each cardinal $\k \geq \o$ the lattice $\la [\k ]^{<\o} ,\subset \ra$ of finite subsets of $\k$ is not reversible.
More generally, if $0\in A \in [\o ]^\o$, then the suborder $\la \bigcup _{n\in A}[\k ]^n ,\subset \ra$ of $[\k ]^{<\o}$ is not reversible.
\end{prop}
\dok
We show that the poset $\P :=\la [\k ]^{<\o} ,\subset \ra$ satisfies the assumptions of Theorem \ref{T8100}, for $\l =\o$.
It is easy to see that $\P$ is well founded of height $\o$, where $L_n =[\k ]^n$, for $n\in \o$, and that $\emptyset$ is the root of $\P$.
$\P$ is $<\o$-directed because the set $[\k ]^{<\o}$ is closed under finite unions.
Also $|P|=|L_1|=|[\k ]^1|=\k$ and $|(\cdot ,p]|<\o$, for all $p\in P$. For the poset $\bigcup _{n\in A}[\k ]^n$ we have a similar proof.
\hfill $\Box$
\begin{ex}\label{EX8100}\rm
By Proposition \ref{T8213} the lattice $\P:=\la [\aleph _\o ]^{<\o} ,\subset \ra$ is not reversible
and, regarding (\ref{EQ8158}) in the proof of Theorem \ref{T8100},
this is witnessed by the b.f.s.\ $\Pi :=\{ \f \in \Pc _{\f ^*}^{<\aleph _\o }(\P ) :  (\dom \f ) \!\downarrow \;= \dom \f \}$.
The desired condensation $f\in \Cond (\P)\setminus \Aut (\P )$ is constructed in Theorem \ref{T8211}
and $\f^* \subset f=\bigcup _{\xi <\aleph _\o}\p _\xi$, where $\p _\xi\in \Pi$, for all $\xi <\aleph _\o$ and, say, $\dom (\f ^*)\subset \o$.
Clearly $f_n:=f\upharpoonright [\o _n]^{<\o} \in \Pi$, for $n\in \o$, but $f=\bigcup _{n\in \o}f_n\not\in \Pi$, because $|f|=\aleph _\o$.
Thus, $\Pi$ is not $|P|$-closed; in fact $\Pi$ is not even $\o _1$-closed ($\s$-closed).
\end{ex}
Let $\N =\{ 1,2,\dots\}$ be the set of natural numbers
and $\mid$ the divisibility relation on $\N$: $m\mid n$ iff $n=km$, for some $k\in \N$.
It is well known that the partial order $\la \N , \,|\, \ra$ is a distributive lattice.
\begin{prop}\label{T8220}
The divisibility lattice $\la \N , \,\mid \, \ra$ is not reversible.
\end{prop}
\dok
We use Theorem \ref{T8100}, for $\k=\l=\o$.
Since $|(\cdot ,n]|<\o$, for all $n\in \N$, the poset $\P:=\la \N , \,\mid \, \ra$ is well founded; since $\P$ is a lattice it is directed; its root is 1.
$L_1$ is the set of all prime numbers
and, generally, for $n\geq 1$, $L_n$ is the set of products of exactly $n$ prime numbers (some of them possibly equal).
Thus $|P|=|L_1|=\he (\P)=\o$.
\hfill $\Box$
\begin{prop}\label{T8221}
For each cardinal $\k \geq \o$ the poset $\Fn (\k,\o)$ of finite partial functions from $\k$ to $\o$ with the partial order $\leq$ given by
$$
f\leq g \;\;\Leftrightarrow \;\;\dom (f) \subset \dom(g) \land \forall \a \in \dom (f) \;\;f(\a )\leq g(\a)
$$
is not reversible.
\end{prop}
\dok
We show that the poset $\P :=\la \Fn (\k,\o) ,\leq \ra$ satisfies the assumptions of Theorem \ref{T8100}, for $\l :=\o$.
Clearly, the empty function, $\emptyset$, is the root of $\P$ and, since $\Fn (\k,\o)\subset [\k \times \o]^{<\o}$ we have $|\Fn (\k,\o)|=\k$.
Also we have $|(\cdot ,f|<\o$, for all $f\in \Fn (\k,\o)$, which implies that the poset $\P$ is well founded.
It is evident that the immediate successors of a function $f\in \Fn (\k,\o)$ are the functions $f\cup \{ \la \a ,0\ra \}$, where $\a \in \k \setminus \dom (f)$,
and the functions $(f\setminus \{ \la \a,f(\a )\ra\}) \cup \{ \la \a,f(\a )+1\ra\}$, where $\a \in \dom (f)$; thus, there are $\k$-many of them.
In particular, the set $L_1 =\{ \la \a,0\ra :\a \in \k\}$ is of size $k$
and generally, for $n\in \o$ we have $L_n=\{ f\in \Fn (\k,\o): |\dom (f)|+\sum \ran (f)=n \}$; so $\he (\P)=\o$.
Finally $\P$ is $<\o$-directed: if $f_i \in \Fn (\k,\o)$, for $i<n$,
and $m:=\max (\bigcup _{i<n}\ran (f_i))$,
then $f:=\{ \la \a ,m\ra :\a \in \bigcup _{i<n}\dom (f_i)\}$ is an upper bound for the set $\{ f_i: i<n\}$.
\kdok
The poset $\la [\k ]^{<\o}, \subset \ra$ witnesses that
the assumptions of Theorem \ref{T8100} are satisfied when $\l=\o$ and $\k$ is any infinite cardinal.
The following example shows that for some $\l >\o$ and $\k \geq\l$ these assumptions can be contradictory in some models of ZFC.
We recall that $\cf ([\k ]^{<\l}):= \min \{ |\CS |: \CS \subset [\k ]^{<\l} \land \forall A\in [\k ]^{<\l} \;\exists S\in \CS \; A\subset S \}$ (see \cite{Bart}, p.\ 14).
\begin{ex}\label{EX8102}\rm
Suppose that $\P$ is a poset satisfying the assumptions of Theorem \ref{T8100}, where $\k =\aleph _\o$ and $\l=\o _1$.
Let $P=\{ p_\a :\a < \aleph _\o\}$ be an enumeration; then $S_\a :=L_1 \cap (\cdot ,p_\a]\in [L_1]^{\leq \o}$, for $\a < \aleph _\o$.
Since $\P$ is $<\o _1$-directed,
for each set $S\in [L_1]^{\leq \o}$ there is $\a < \aleph _\o$ such that $S \subset (\cdot ,p_\a]$
and, hence, $S\subset S_\a$, that is, $S\in P(S _\a)$.
So $\cf ([\aleph _\o]^{\leq \o}])=\aleph _\o$, we have $[L_1]^{\leq \o}=\bigcup _{\a < \aleph _\o}P(S _\a)$
and, hence,  $\aleph _\o ^\o =|[L_1]^{\leq \o}|\leq \sum _{\a < \aleph _\o}|P(S _\a)|\leq \aleph _\o \fc$.
Thus, if $\fc <\aleph _\o$, then $\aleph _\o ^\o \leq \aleph _\o $, which is impossible by K\"{o}nig's theorem.
Moreover $\cf ([\aleph _\o]^{\leq \o}])\fc =\aleph _\o ^\o > \aleph _\o$ (see \cite{AbrMag} p.\ 1196); so if $\fc <\aleph _\o ^\o$ we have a contradiction again.
\end{ex}
In the sequel we present more applications of Theorem \ref{T8100}.
We recall that if $\vartriangleleft$ is a binary relation on a set $P$ without loops
(for each $n\in \N$ and $p_0, \dots ,p_n \in P$ such that $p_0 \vartriangleleft p_1 \vartriangleleft \dots \vartriangleleft p_n$ we have $p_0\neq p_n$),
then the transitive closure, $<$, of $\vartriangleleft$ is given by: $p < q$ iff there are $n\in \N$ and $p_0,\dots,p_n \in P$ such that
$p=p_0 \vartriangleleft p_1 \vartriangleleft \dots \vartriangleleft p_n=q$.
Then $<$ is a strict partial order on the set $P$.
Let $\o \leq\l \leq \k$ be cardinals, $P:=\bigcup _{\a <\l}P_\a$, where

$P_0 =\{ 0\}\times\k =\{\la 0,\xi\ra :\xi \in \k \}$,

$P_\b =\{ \b \}\times ([\bigcup _{\a <\b}P_\a]^{< \l}\setminus \{ \emptyset \}) $, for $0<\b <\l$,

\noindent
let $\vartriangleleft$ be the binary relation on $P$ given by $\la \a ,S\ra \vartriangleleft \la \b ,T\ra$ iff $\a < \b$ and $\la \a ,S\ra  \in T$,
and let $<$ be the transitive closure of $\vartriangleleft$, namely,
$\la \a ,S\ra < \la \b ,T\ra$ iff there are $n\in \N$ and $\la \a _0 ,S_0\ra, \dots ,\la \a _n ,S_n\ra \in P $ such that

(i) $\a =\a _0 <\dots <\a_{n-1} < \a _n =\b$ and

(ii) $\la \a ,S\ra = \la \a_0 ,S_0\ra \in S_1 \land \dots \land \la\a_{n-1} ,S_{n-1}\ra \in S_n \land \la\a_n ,S_n\ra=\la \b ,T\ra$.

\noindent
Then $\P _{\l,\k}:=\la P,<\ra$ is a strict partial order.
\begin{prop}\label{T8218}
(GCH) Let $\l \leq \cf (\k )<\k$, where $\l$ is a regular cardinal.
Then the poset $\P _{\l,\k}$ satisfies the assumptions of  Theorem \ref{T8100} and it is not reversible.
\end{prop}
\dok
By the assumptions, $|[\k]^{<\l}| =\sum _{\mu <\l} |[\k ]^\mu|=\sum _{\mu <\l} \k ^\mu =\l \k =\k$.

Assuming that $\la \a _0 ,S_0\ra >\la \a _1 ,S_1\ra > \dots$ is a decreasing chain in $\P _{\l,\k}$
by (i) we would have $\a _0 >\a _1 > \dots$, which is impossible. So $\P _{\l,\k}$ is well founded.

Since $|P_0|=\k$ we have $|P|\geq\k$.
In order to show that $|P|\leq\k$
by induction we prove that $|P_\a|\leq \k$, for all $\a <\l$.
Let $\b <\l$ and $|P_\a|\leq \k$, for all $\a <\b$.
Then $|\bigcup _{\a <\b}P_\a|=\sum _{\a <\b}|P_\a|\leq =\k$;
so, $|P_\b| =|[\bigcup _{\a <\b}P_\a]^{<\l}|\leq \k^{<\l }=\k$.

If $\la 0,\xi \ra \in P_0$, then by (i) $\la \a ,S\ra <\la 0,\xi \ra$ would give $\a <0$, which is false.
Thus $P_0 \subset\Min (\P _{\l,\k})$ and, hence, $|\Min (\P _{\l,\k})|=\k$.

$\P _{\l,\k}$ is $< \l$-directed.
Let $\mu <\l$ and $\la \a _\xi, S_\xi \ra\in P_{\a _\xi}$, for $\xi <\mu$,
and let $\b <\l$, where $\a _\xi <\b$, for all $\xi <\mu$.
Then $T:=\{ \la \a _\xi, S_\xi\ra :\xi <\mu\}\in [\bigcup _{\a <\b}P_\a]^{<\l}\setminus \{ \emptyset \}$
and, hence, $\la \b,T\ra \in P_\b$.
For $\xi <\mu$ we have $\a _\xi <\b$ and $\la \a _\xi, S_\xi \ra \in T$ so $\la \a _\xi, S_\xi \ra \vartriangleleft \la \b,T\ra $ and, hence, $\la \a _\xi, S_\xi \ra <\la \b,T\ra $.

By induction we prove that for each $\b <\l$ and $p\in P_\b$ we have $|(\cdot ,p]|<\l $.
For $\b =0$ this is true, because $P_0\subset \Min (\P _{\l,\k})$.
Let $\b <\l$ and suppose that
\begin{equation}\label{EQ8179}
\forall \a <\b \;\;\forall p\in P_\a \;\;|(\cdot ,p]|<\l.
\end{equation}
Let $p=\la \b,T\ra \in P_\b$.
Then $T=\{ p_i :i \in I\}\subset \bigcup _{\a <\b}P_\a$, $0<|I| <\l$,
and for each $i\in I$ we have $p_i =\la \a _i ,S_i\ra$, where $\a _i <\b$,
and, by (\ref{EQ8179}), $|(\cdot ,p_i]|<\l$.
Thus, it remains to be proved that $(\cdot ,p)\subset \bigcup _{i\in I}(\cdot ,p_i]$.
So, if $\la \a,S\ra < \la \b,T\ra$,
then there are $n\in \N$ and $\la \a _0 ,S_0\ra, \dots ,\la \a _n ,S_n\ra \in P $ such that

(i) $\a =\a _0 <\dots <\a_{n-1} < \a _n =\b$ and

(ii) $\la \a ,S\ra = \la \a_0 ,S_0\ra \in S_1 \land \dots \land \la\a_{n-1} ,S_{n-1}\ra \in S_n =T$,

\noindent
By (ii) we have  $\la\a_{n-1} ,S_{n-1}\ra =p_i =\la \a _i ,S_i\ra$, for some $i\in I$.
So, by (i) and (ii) if $n\geq 2$ we have $\a =\a _0 <\dots  < \a_{n-1}=\a_i $ and
$\la \a ,S\ra = \la \a_0 ,S_0\ra \in S_1, \dots , \la\a_{n-2} ,S_{n-2}\ra \in S_{n-1},  \la\a_{n-1} ,S_{n-1}\ra=\la \a _i ,S_i\ra $,
which means that $\la \a,S\ra \in (\cdot , p_i )$.
For $n= 1$ we have  $\a =\a _0 < \a _1 =\b$ and $\la \a ,S\ra = \la \a_0 ,S_0\ra \in S_1 =T$,
thus $\la \a ,S\ra\in T$,
which gives $\la \a ,S\ra=p_i$, for some $i\in I$ and, again,  $\la \a,S\ra \in (\cdot , p_i ]$.
So $(\cdot ,p]=\bigcup _{i\in I}(\cdot ,p_i]$.

Finally we show that $\he (\P _{\l,\k} )= \l$.
Let us define $p_0 =\la 0,0\ra$ and $p_\b = \la \b ,\{ p_\a :\a <\b\}\ra$, for $0<\b<\l$.
Then $\la p_\b :\b <\l\ra$ is a $\vartriangleleft$-chain of type $\l$
and, hence $\he (\P _{\l,\k} )\geq \l$.
Assuming that $\he (p)=\l$, for some $p\in P$
by Fact \ref{T8223}(b) we would have $|(\cdot ,p]|\geq \l$, which is false.
\hfill $\Box$
\section{Approximations with open  domain and bounded field. More ideals}\label{S6}
If $\P =\la P,\leq\ra$ is a poset and $p,q\in P$, then $q$ is an {\it immediate successor of $p$}, we write $q\in \Is (p)$, iff $p<q$ and there is no $r\in P$ such that $p<r<q$
\begin{te}\label{T8206}
If $\P$ is a  $<\!\k$-directed poset of (regular) size $\k\geq \o$, with a root $r$, $\Is (r) \setminus (\cdot ,p]\neq \emptyset$, for all $p\in P$, some $a_0\in \Is (r)$ has an immediate successor and
\begin{equation}\label{EQ8166}
\forall a,p\in P \;\;\exists \eta \in \Pc (\P ) \;\;\eta :(\cdot,a] \rightarrow [p,\cdot ),
\end{equation}
then $\P$ is not reversible.
\end{te}
\dok
Let $b_0 \in \Is (a_0)$ and $a_1\in \Is (r)\setminus \{ a_0\}$.
Clearly the mapping $\f ^* =\{ \la r ,r \ra,\la a_0, b_0 \ra,\la a_1, a_0 \ra\} \in \Pc (\P )$ is a bad condensation,
$\dom \f ^* =\{ r, a_0 , a_1\}$ is an open set, $\ran \f ^* =\{ r, a_0 , b_0\}$
and, since $\P$ is a directed poset, the set $\dom \f \cup \ran \f$ has an upper bound.
We note that by (\ref{EQ8166}) a set having an upper bound has a strict upper bound.
We show that $\Pi$ is a b.f.s.\ where
\begin{equation}\label{EQ8164}\textstyle
\Pi :=\{ \f \in \Pc _{\f ^*}(\P ) : (\dom \f ) \!\downarrow = \dom \f \land \exists p\in P \;\; \dom \f \cup \ran \f \subset (\cdot , p] \}.
\end{equation}

(bf1) Let $\f \in \Pi$ and $a\in P\setminus \dom \f$.
It is evident that the set $D=\dom \f \cup (\cdot ,a ]$ is open.
By (\ref{EQ8164}) there is $p\in P$ such that $\dom \f \cup\ran \f < p$,
and by (\ref{EQ8166}) there is $\eta \in \Pc (\P )$, where $\eta :(\cdot,a] \rightarrow [p,\cdot )$.
Then $q:= \eta (a)> p$ and $\eta [(\cdot,a]]\subset [p,q]$.
Clearly the mapping $\p :D\rightarrow P$ defined by $\p =\f \cup \eta \upharpoonright ((\cdot ,a]\setminus \dom \f)$ is an injection.

We prove that $\p$ is a homomorphism; let $x,y\in D$ and $x<y$.
Since the set $\dom \f$ is open it is impossible that $x\not\in \dom \f \ni y$.
If $x,y\in \dom \f $, then, since $\f$ is a homomorphism, we have $\p (x)=\f (x)< \f (y)=\p (y)$.
If $x,y\not\in \dom \f $, then $\p (x)=\eta (x)< \eta (y)=\p (y)$, because $\eta$ is a homomorphism.
If $x\in \dom \f \not\ni y$, then $y\in (\cdot ,a]\setminus \dom \f$
and, hence, $\p (x)=\f (x)<p \leq \eta (y)=\p (y)$.

Thus $\p \in \Pc (\P )$,
$\f ^* \subset \f \subset \p$
and the set $\dom \p =\dom \f \cup (\cdot ,a ]$ is open.
Since $\P$ is directed there is $p'\in P$ such that $a,q < p'$,
and hence, $\dom \p \cup \ran \p < p'$.
So $\p \in \Pi$.
Since $\f \subset \p$ and $a\in \dom \p$, (bf1) is true.

(bf2) Let $\f \in \Pi$ and $b\in P\setminus \ran \f$.
By (\ref{EQ8164}) there is $p\in P$ such that $\dom \f \cup\ran \f < p$.
Since $\Is (r) \setminus (\cdot ,p]\neq \emptyset$ we can take $a\in \Is (r) \setminus (\cdot ,p]$
So $a\not\in \dom \f$ and $\p := \f \cup \{ \la a ,b \ra\}$ is an injection.
Since $r\in \dom \f$ the set $\dom \p =\{ a\} \cup \dom \f$ is open.

We show that $\p$ is a homomorphism.
Let $x,y\in \{ a\} \cup \dom \f$ and $x < y$.
Since the set $\dom \f$ is open it is impossible that $x\not\in \dom \f \ni y$;
since $\dom \p \setminus \dom \f =\{ a\}$, $x,y\not\in \dom \f $ is impossible too.
If $x,y\in \dom \f $, then, again, $\p (x)=\f (x)< \f (y)=\p (y)$.
If $x\in \dom \f \not\ni y$,
then $y=a$
and, hence $x=r$.
Now, since $b \not\in \ran \f$
and, hence, $b\neq r$,
we have $\p (x)=\f ^*(r)= r< b=\p (a)=\p (y)$.

Thus $\p \in \Pc (\P )$, $\f ^* \subset \f \subset \p$
and the set $\dom \p =\dom \f \cup \{ a \}$ is open.
Since the poset $\P$ is directed there is $p'\in P$ such that $a,b,p <p'$, and hence, $\dom \p \cup \ran \p < p'$.
So $\p \in \Pi$.
Since $\f \subset \p$ and $b\in \ran \p$, (bf2) is true. So, $\Pi$ is a b.f.s.\ indeed.

(b) By Theorem \ref{T8196} it remains to be proved that the poset $\la \Pi ,\supset\ra$ is $\k$-closed.
So let $\g \in \k$ and let $\f _\a \in \Pi$, where $\f _0 \subset \f _1 \subset \dots$.
Then, clearly,  $\f ^* \subset \f :=\bigcup _{\a <\g}\f _\a \in \Pc (\P)$
and the set $\dom \f :=\bigcup _{\a <\g}\dom \f _\a$ is open.
By (\ref{EQ8164}) for each $\a <\g$ there is $p_\a \in P$ such that $\dom \f _\a \cup\ran \f _\a < p_\a$
and, since $\P$ is $<\k$-directed and $|\g|<\k$,
there is $p \in P$ such that $p_\a \leq p$, for all $\a <\g$.
But then $\dom \f  \cup\ran \f  =\bigcup _{\a <\g}\dom \f _\a \cup\ran \f _\a < p$;
thus $\f \in \Pi$ and $\f$ is a lower bound for the sequence $\la \f _\a :\a <\g\ra$.
\hfill $\Box$
\begin{rem}\label{R8100}\rm
Concerning Theorem \ref{T8206} we note that

(a) The theorem is true if $\P$ has no root and if we replace $\Is (r)$ by $\Min (\P )$;
In addition, these two variants of Theorem \ref{T8100} have their obvious dual versions.

(b) Condition (\ref{EQ8166}) is satisfied if the poset $\P$ is homogeneous in the following sense: $[p ,\cdot) \cong \P$, for each $p\in P$;

(c) The assumptions of the theorem imply that $|\Is (r)|=\k$ (otherwise there would be $p>\Is (r)$ and we would have $\Is (r) \setminus (\cdot ,p]= \emptyset$)
and that $\k$ is a regular cardinal (otherwise, there would be $p:=\max P$ and $\Is (r) \setminus (\cdot ,p]= \emptyset$).
\end{rem}
\begin{ex}\label{EX8103}\rm
The suborder $\P =\la P,\leq \ra$ of the integer plane, $\Z ^2$,
where $P=\{ \la m,n\ra \in Z^2: m+n\geq 0\}$,
satisfies the assumptions of Theorem \ref{T8206} (for $\k=\o$); so, $\P$ is not reversible.
($\P$ is directed, has no root, $\Min (\P )=\{ \la m,-m\ra : m\in Z\}$, etc.)
\end{ex}
\begin{prop}\label{T8121}
If $\o \leq \l =\l ^{< \l}\leq \k$, then the poset $\la [\k ]^{<\l }, \subset \ra$ is not reversible.
Under GCH, the poset $\la [\k ]^{<\l }, \subset \ra$ is not reversible, whenever $\o \leq \l \leq \k$ and $\l$ is a regular cardinal.
\end{prop}
\dok
We recall that an infinite cardinal $\k$ such that $\k ^{<\k}=\k$ is regular
($\cf  (\k )<\k$ would imply that $\k ^{\cf  (\k )}\leq \k$, which is false by K\"{o}nig's theorem).
Under the GCH these two conditions are equivalent; on the other hand, $\o _1 {}^{<\o _1}>\o_1$ in each model of $\neg$ CH.

First, let $\l=\k$. We show that the poset $\P =\la [\k ]^{<\k } ,\subset \ra$ satisfies the conditions of Theorem \ref{T8206}.
By the assumptions,  the poset $\P $ is of size $\k$
and $\k$ is a regular cardinal.
So, the union of $<\k$-many $<\k$-sized subsets of $\k$ is of size $<\k$ and, hence, $\P$ is $<\k$-directed.
The root $\emptyset$ of $\P$ has $\k$-many immediate successors: $\{ a\}$, $\a \in \k$, and each $p\in P$ covers $<\k$ of them.
If $a,p\in P$, then $a\in [\k ]^{<\l}$ and $p\in [\k ]^{<\mu}$, for some cardinals $\l ,\mu <\k$,
and $(\cdot ,a]=P(a)$. Let $b\in [\k ]^{<\l}$, where $b\cap p=\emptyset$, and let $q=p\cup b$.
Then $q\in P$, $[p,q]\cong P(b)\cong P(a)$ and (\ref{EQ8166}) is true.

Second, let $\l < \k$.
Then, as above, there is $f\in \Cond (\la [\l ]^{<\l }, \subset \ra) \setminus \Aut (\la [\l ]^{<\l }, \subset \ra)$.
If $A\in [\k ]^{<\l }$,
then $A\cap \l\in [\l ]^{<\l }$;
so, the function $F: [\k ]^{<\l }\rightarrow [\k ]^{<\l }$ given by
\begin{equation}\label{EQ8167}
F(A)=f(A\cap \l)\cup (A\setminus \l), \;\;\mbox{ for } A\in [\k ]^{<\l },
\end{equation}
is well defined.
{\it $F$ is an injection.}
If $F(A)=F(B)$,
then $F(A)\cap \l =F(B)\cap \l $
and, by (\ref{EQ8167}), $f(A\cap \l)= f(B\cap \l)$,
which, since $f$ is an injection, implies  $A\cap \l= B\cap \l$.
Since in addition we have $F(A)\setminus \l =F(B)\setminus \l $,
that is, by (\ref{EQ8167}), $A\setminus \l =B\setminus \l $,
we finally have $A=B$.
{\it $F$ is an surjection.}
If $B\in [\k ]^{<\l }$, then $B\cap \l \in [\l ]^{<\l }$
and, since $f$ is onto,
there is $A'\in [\l ]^{<\l }$ such that $B\cap \l=f(A')$.
Now $A:= A' \cup (B\setminus \l)\in [\k ]^{<\l }$
and $F(A)=f(A')\cup (B\setminus \l)=B$.
{\it $F$ is a homomorphism.}
If $A, B\in [\k ]^{<\l }$ and $A\subset B$,
then $A\cap \l\subset B\cap \l$
and, since $f$ is a homomorphism, $f(A\cap \l)\subset f(B\cap \l)$.
In addition we have $A\setminus \l\subset B\setminus \l$
and, by (\ref{EQ8167}), $F(A)\subset F(B)$.
So, $F\in \Cond (\la [\k ]^{<\l }, \subset \ra) $
and, since $F\upharpoonright [\l ]^{<\l }=f$, we have $F\not\in \Aut (\la [\k ]^{<\l }, \subset \ra) $.
Thus, the poset $\la [\k ]^{<\l }, \subset \ra$ is not reversible indeed.
\kdok
Further we consider the direct power $\la {}^\k \k ,\leq \ra$, where for $f,g:\k \rightarrow \k $ we have: $f\leq g$ iff $f(\a )\leq g(\a)$, for all $\a <\k$.
For $f:\k \rightarrow \k $ let $\supp (f):=\{ \a \in \k : f(\a )>0\}$.
\begin{prop}\label{T8222}
Let $\o \leq \k =\k ^{< \k}$ and $P:=\{ f\in {}^\k \k : |\supp (f)|<\k \}$. Then the suborder $\P :=\la P, \leq \ra$ of $\la {}^\k \k ,\leq \ra$ is not reversible.
\end{prop}
\dok
We show that the poset $\P$ satisfies the conditions of Theorem \ref{T8206}.
Since $|[\k]^{<\k}|=\k ^{< \k}=\k$ there are $\k$-many different supports and if $S\in [\k ]^\l$ is one of them, where $\l <\k$,
there are $\k ^\l =\k$ different functions $f\in P$ such that $\supp (f)=S$; thus $|P|=\k$,
in addition, by the regularity of $\k$, the set $\ran (f)$ is bounded in $\k$.
If $\mu <\k$ and $f_\xi \in P$, for $\xi <\mu$,  then,  by the regularity of $\k$ again, the set $D :=\bigcup _{\xi <\mu}\supp (f_\xi )$ is of size $\k$
and there is $\b :=\sup (\bigcup _{\xi < \mu} \ran (f_\xi))$;
so, defining $f \in P$ by $f(\a)=\b$, for $\a \in D$, and $f(\a)=0$, for $\a \in \k \setminus D$,
we obtain an upper bound for the set $\{ f_\xi :\xi <\mu \}$; thus $\P$ is $<\k$-directed.

Clearly, the ``zero function" $f_0$ is the root of $\P$ and $\Is (f_0)=\{ g_\xi :\xi <\k\}$, where $g_\xi (\xi )=1$ and $g(\a )=0$, for $\a \neq \xi$.
Defining $h(0)=2$ and $h(\a )=0$, for $\a \neq 0$, we have $h\in \Is (g_0)$
For $f\in P$ we have $|(\cdot ,f] \cap \Is (f_0)| =|\{ g_\xi :\xi\in \supp (f)\}|<\k$;
thus, since  $|\Is (f_0)|=\k$ we have $|\Is (f_0)\setminus (\cdot ,f]|=\k$.

In order to prove (\ref{EQ8166}) we show that $[f ,\cdot ) \cong \P$, for arbitrary $f\in P$.
For $f,g\in P$ let the function $f+g \in P$ be defined by $(f+g)(\a)=f(\a) +g(\a)$, for all $\a <\k$.
We show that the mapping $F:P\rightarrow [f ,\cdot )$ defined by $F(g)=f+g$ is an isomorphism.

{\it $F$ is an injection.} If $F(g)=F(h)$, then for each $\a <\k$ we have $f(\a) +g(\a)=f(\a) +h(\a)$
and, since the ordinal addition is left-cancelative, $g(\a)=h(\a)$; thus $g=h$.

{\it $F$ is a surjection.} If $h\in [f ,\cdot )$, then for each $\a <\k$ we have $f(\a) \leq h(\a)$
and there is a unique ordinal $\b _\a $ such that $f(\a) +\b _\a = h(\a)$.
Clearly, for $\a \in \k \setminus (\supp (f)\cup \supp (h))$ we have $\b _\a=0$
and defining $g(\a )=\b _\a $, for all $\a <\k$, we have $\supp (g)\subset \supp (f)\cup \supp (h)$
we have $g\in P$.
In addition $F(g)(\a )=f(\a )+\b _\a = h(\a)$, for all $\a <\k$, and, hence, $F(g)=h$.

{\it $F$ is an isomorphism.} For $g,h\in P$ we have $F(g)\leq F(h)$ iff for each $\a <\k$ we have $f(\a) +g(\a)\leq f(\a) +h(\a)$,
which is (by the left cancelation for inequalities) equivalent to $g(\a)\leq h(\a)$, iff $g\leq h$.
\kdok
Now we consider the $\s$-ideal of meager (first category) sets in the Borel algebra of a Polish space.
If $\X =\la X,\CO\ra$ is a topological space, $\Borel (\X )$ will denote the Boolean algebra of Borel subsets of $X$ and  $\CM (\X )$  the ideal of meager subsets of $X$.
We will use the following fact.
\begin{fac}\label{T8207}
Let $\X =\la X,\CO\ra$ be a topological space, $\emptyset \neq A\subset X$ and $\A =\la A, \CO _A\ra$ the corresponding subspace. Then

(a) $\Nwd (\A )\subset \{ N\cap A : N\in \Nwd (\X )\}\subset \Nwd (\X )$ and $\CM (\A )\subset \CM (\X )$;

(b) If $A\in \Borel (\X )$, then $\Borel (\A )=\{ B\cap A: B\in \Borel (\X )\}\subset \Borel (\X )$;

(c) If $h: \X \rightarrow \A$ is a homeomorphism, then $\Borel (\A )=\{ h[B]:B\in \Borel (\X )\}$ and $\CM (\A )=\{ h[M]:M\in \CM (\X )\}$;

(d) If $\X$ is a Baire space, each meager set is contained in a Borel meager set.
\end{fac}
\dok
(a) For $N\in \Nwd (\A )$ we have $N=N\cap A$.
Suppose that $N\not \in \Nwd (\X )$ and let $x\in O \in \CO$, where $O\subset \overline{N}^\X $.
Then $O\in \CU_\X (x)$ and, hence, $O\cap N\neq \emptyset$.
Thus $\CO _\A \setminus \{ \emptyset \}\ni O\cap A \subset \overline{N}^\X \cap A =\overline{N}^\A$
which gives $N\not \in \Nwd (\A )$ and we have a contradiction. Clearly, if $N\in \Nwd (\X )$, then $N\cap A \in \Nwd (\X )$.

If $M\in \CM (\A )$, then $M=\bigcup _{n\in \o}N_n$, where $N_n\in \Nwd (\A )$, and by (a) we have $N_n\in \Nwd (\X )$, for all $n\in \o$. Thus $M\in \CM (\X )$.

(b) Let $A\in \Borel (\X )$. Since $\CO \subset \Borel (\X )$ and the family $\Borel (\X )$ is closed under $\cap$, we have $\CO _A =\{ O\cap A : O\in \CO\}\subset \Borel (\X )$.
$\Borel (\A )$ is the minimal closure of $\CO _A$ under countable intersections and complements. So, since the family $\Borel (\X )$ is closed under countable intersections and complements,
we have $\Borel (\A )\subset \Borel (\X )$.

(c) For the first equality see \cite{Eng} p.\ 36. The second is evident.

(d) Let $M=\bigcup _{n\in \o}N_n$, where $\Int (\overline{N_n})=\emptyset$, for all $n\in \o$. Then $M\subset \bigcup _{n\in \o}\overline{N_n}\in\Borel (\X )\cap \CM (\X )$.
\kdok
We recall that the cardinal invariant $\add (\CM (\X)):=\min\{ |\CM|: \CM \subset \CM (\X) \land \bigcup \CM \not \in \CM (\X)\}$
is the same for several ``spaces of reals", e.g.\ $\BR$, $2^\o$, $\o ^\o$ (see \cite{Bart}, p.\ 16), it is usually denoted by $\add (\CM )$.
The equality $\add (\CM )=\mathfrak{c}$ follows from CH, MA and PFA.
\begin{prop}\label{T8208}
($\add (\CM )=\mathfrak{c}$) The poset of meager Borel subsets of the Baire space $\o ^\o$ is not reversible.
\end{prop}
\dok
We show that this poset satisfies the assumptions of Theorem \ref{T8206}.
Let $\add (\CM )=\mathfrak{c}$ hold, $X:=\o ^\o$,
let $\X :=\la X, \CO \ra$ be the Baire space and $P=\Borel (\X )\cap \CM (\X )$.
Then ${\mathfrak c}$ is a regular cardinal (see \cite{Bart}, p.\ 12),
and since $[X]^1 \subset P$ and $|\Borel (\X )|={\mathfrak c}$, the poset $\P =\la P,\subset\ra$ is of size ${\mathfrak c}$.
If $\l <{\mathfrak c}$ and $p_\a \in P$, for $\a <\l$,
then since  $\l<\add (\CM )$ we have $M=\bigcup _{\a <\l}p_\a \in \CM (\X )$
and by Fact \ref{T8207}(d) there is $p\in P$ such that $M\subset p$;
thus, the poset $\P$ is $(< \! {\mathfrak c})$-directed.
The empty set is the root of $\P$ and $\Is (\emptyset)=\{ \{ x \} :x\in X\} \in [P]^{\mathfrak c}$.
If $p\in P$, then $\Is (\emptyset) \subset (\cdot ,p]$ would imply that $p=X$, which is false by the Baire theorem.

Let $a,p\in P$.
Then $p=\bigcup _{n\in \o} r_n$, where $\Int _\X ( \overline{r_n})=\emptyset$
and $p\subset p':=\bigcup _{n\in \o} \overline{r_n}$, where $\overline{r_n}$, $n\in \o$, are closed nowhere dense sets.
Thus the set $G := X\setminus p'$ is a dense co-meager $G_\delta$-set in $\X$.
Also, there is a Borel set $A\subset G$ homeomorphic to $\X$; namely, if the set $p'$ is dense, then $G\cong \X$ (see \cite{Kech}, p.\ 40) and we take $A=G$;
otherwise, there is a basic clopen set $A\subset G$ and, clearly $A\cong \X$.
Then by (a) and (b) of Fact \ref{T8207} we have  $\CM (\A )\subset \CM (\X )$ and $\Borel (\A )\subset \Borel (\X)$; thus
\begin{equation}\label{EQ8168}
\Borel (\A )\cap \CM (\A ) \subset \Borel (\X )\cap \CM (\X ).
\end{equation}
Let $h: \X \rightarrow \A$ be homeomorphism
If $b\in (\cdot ,a]$, then $b\in \Borel (\X )\cap \CM (\X )$,
by (c) of Fact \ref{T8207} we have $h[b]\in \Borel (\A )\cap \CM (\A )$,
and, by (\ref{EQ8168}), $h[b]\in \Borel (\X )\cap \CM (\X )$.
Thus, since $p\in \Borel (\X )\cap \CM (\X )$, we have
\begin{equation}\label{EQ8169}
\forall b\in (\cdot ,a]\;\; p\cup h[b]\in \Borel (\X )\cap \CM (\X )
\end{equation}
and, in particular, $q:=p \cup h[a]\in \Borel (\X )\cap \CM (\X )$.
We show that the mapping $f:(\cdot ,a]\rightarrow [p, q]$ defined by
$$
f(b)=p\cup h[b] , \;\; \mbox{ for all }b\in (\cdot ,a],
$$
is an isomorphism.
For $b\in (\cdot ,a]$ by (\ref{EQ8169}) we have $f(b)\in \Borel (\X )\cap \CM (\X )=P$
and, since $b\subset a$ we have $p\subset p\cup h[b] =f(b)\subset p\cup h[a]=q$;
thus, $f$ maps $(\cdot ,a]$ to $[p, q]$ indeed.

{\it $f$ is an injection.} If $b,b'\in (\cdot ,a]$ and $b\neq b'$, then since $h$ is a bijection we have $h[b]\neq h[b']$ and, hence, $f(b)\neq f(b')$.

{\it $f$ is a surjection.} If $r\in [p,q]$,
then $r\in \Borel (\X )\cap \CM (\X )$
and $p\subset r \subset p \cup h[a]$.
Thus $r\setminus p \subset h[a]\subset A$,
and, since $r\setminus p\in \Borel (\X )$,
by Fact \ref{T8207}(b) we have $r\setminus p = (r\setminus p)\cap A\in \Borel (\A )$.
So, by Fact \ref{T8207}(c) there is $b\in \Borel (\X )$ such that $h[b]=r\setminus p$,
and, since $r\setminus p \subset h[a]$, we have $b\subset a$.
Since $a\in \CM (\X )$ by Fact \ref{T8207}(c) we have $h[a]\in \CM (\A )$
and $r\setminus p \subset h[a]$ gives $r\setminus p\in \CM (\A )$.
So, by Fact \ref{T8207}(c) there is $b'\in \CM (\X )$ such that $h[b']=r\setminus p$.
Now $h[b]=h[b']$ gives $b=b' \in \Borel (\X )\cap \CM (\X )$.
Thus $b\in (\cdot ,a]$ and $r=p\cup (r\setminus p)= p\cup h[b]=f(b)$.

{\it $f$ is a isomorphism.} If $b,b'\in (\cdot ,a]$, then since $h$ is a bijection we have: $b\subset b'$ iff $h[b]\subset h[b']$ iff $f(b)\subset f(b')$.
\kdok
The results concerning the ideals in Boolean algebras can be converted into the corresponding results about their dual filters. For example, by Proposition \ref{T8213} for each cardinal $\k \geq \o$
the Fr\'{e}chet filter of cofinite subsets of $\k$ is not reversible.
\section{Approximations with convex  domain and bounded field. Direct powers of chains}\label{S7}
If $\k$ is a cardinal, a poset $\P$ will be called {\it downwards $<\k$-directed} iff each set $S\subset P$ of size $<\k$ has a lower bound.
\begin{te}\label{T8209}
If $\P$ is a poset of (regular) size $\k \geq \o$, which is both upwards and downwards $<\k$-directed and if
\begin{equation}\label{EQ8174}
\forall p,q\in P \;\;(p\leq q \Rightarrow \exists a \in P \; a \parallel [p,q]),
\end{equation}
\begin{equation}\label{EQ8175}
\forall p,q,r\in P\;\; \exists \eta \in \Pc (\P) \;\;\eta :[p,q] \rightarrow (\cdot ,r],
\end{equation}
\begin{equation}\label{EQ8176}
\forall p,q,r\in P\;\; \exists \eta \in \Pc (\P) \;\;\eta :[p,q] \rightarrow [r,\cdot ),
\end{equation}
\noindent
then $\P$ is not reversible.
\end{te}
\dok
By (\ref{EQ8174}) there are incomparable $a_0,a_1\in P$.
Since $\P$ is upwards directed there is $b_0\geq a_0,a_1$ and, hence, $b_0 >a_0$.
Clearly, the mapping $\f ^* :=\{ \la a_0,b_0\ra,\la a_1,a_0\ra\}$ is a bad condensation,
$\dom (\f ^*)=\{ a_0,a_1\}$ is a convex set, $\ran (\f ^*)=\{ a_0,b_0\}$
and, since $\P$ is upwards and downwards directed, there are $p,q\in P$ such that $\dom (\f ^*)\cup\ran (\f ^*)\subset [p,q]$.

By (\ref{EQ8175}), for each $r\in P$ there is $\eta \in \Pc (\P)$ such that $\eta :[a_0,b_0] \rightarrow (\cdot ,r]$
which gives $\eta (a_0)<r$.
Thus $\P$ has no minimal elements
and, by (\ref{EQ8176}), $\P$ has no maximal elements.
Consequently, if a set $S\subset P$ is bounded, then there are $p,q\in P$ such that $S\subset (p,q)$.

We show that $\Pi$ is a b.f.s.\ where
\begin{equation}\label{EQ8170}\textstyle
\Pi :=\{ \f \in \Pc _{\f ^*} (\P ) : \dom \f  \mbox{ is convex } \land \exists p,q\in P \;\; \dom \f \cup \ran \f \subset (p,q) \};
\end{equation}
(bf1) Let $\f \in \Pi$, $a\in P\setminus \dom \f$ and  $\dom \f \cup \ran \f \subset (p,q)$.
Since the set $\dom \f$ is convex it is impossible that both $\dom \f \cap (\cdot ,a )\neq \emptyset$ and $\dom \f \cap (a,\cdot )\neq \emptyset$.

{\it Case 1}: $\dom \f \cap (\cdot ,a )\neq \emptyset$. First we show that the set $D=\dom \f \cup ((\dom \f )\uparrow  \cap (\cdot ,a ])$ is convex.
Let $r,s \in D$ and $r \leq t \leq s$.
If $r,s\in \dom \f$, then $t\in \dom \f \subset D$, because $\dom \f$ is convex.
If $r,s\not\in \dom \f$,
then $s\leq a$
and there is $r'\in \dom \f$ such that $r' \leq r$;
so, $r' \leq r \leq t \leq s\leq a$ and, hence, $t\in D$.
If $r\in \dom \f\not\ni s$,
then $s\leq a$
so, $r \leq t \leq s\leq a$ and, hence, $t\in D$.
If $r\not\in \dom \f \ni s$,
then there is $r'\in \dom \f$ such that $r' \leq r$;
and $r' \leq r\leq s$; but this is impossible, by the convexity of $\dom \f$.

For $r\in D\setminus \dom \f$ we have $r\leq a$
and there is $r'\in \dom \f$ such that $r' \leq r$
and, since $\dom \f \subset (p,q)$, we have $r\in (p,a]$.
Thus $D\setminus \dom \f \subset [p,a]$,
by (\ref{EQ8176}) there is $\eta \in \Pc (\P )$, where $\eta :[p,a] \rightarrow [q,\cdot )$.
Clearly the mapping $\p :D\rightarrow P$ defined by $\p =\f \cup \eta \upharpoonright (D\setminus \dom \f)$ is an injection.

We prove that $\p$ is a homomorphism;
let $x,y\in D$ and $x<y$.
As above, is impossible that $x\not\in \dom \f \ni y$.
If $x,y\in \dom \f $, then, since $\f$ is a homomorphism, we have $\p (x)=\f (x)< \f (y)=\p (y)$.
If $x,y\not\in \dom \f $, then $\p (x)=\eta (x)< \eta (y)=\p (y)$, because $\eta$ is a homomorphism.
If $x\in \dom \f \not\ni y$, then $y\in D\setminus \dom \f$
and, hence, $\p (x)=\f (x)<q \leq \eta (y)=\p (y)$.

Thus $\p \in \Pc (\P )$,
$\f ^* \subset \f \subset \p$
and the set $\dom \p =D$ is convex.
Since $\P$ is upwards directed there is $q'\in P$ such that $a,\eta (a) < q'$,
and hence, $\dom \p \cup \ran \p < q'$.
Also we have $p<\dom \p \cup \ran \p $;
thus $\dom \p \cup \ran \p\subset (p,q')$
and $\p \in \Pi$.
Since $\f \subset \p$ and $a\in \dom \p$, (bf1) is true.

{\it Case 2}: $\dom \f \cap (a,\cdot )\neq \emptyset$.
The proof is dual:
the set $D: =\dom \f \cup ((\dom \f )\!\downarrow  \cap [a,\cdot ))$ is convex,
$D\setminus \dom \f \subset [a,q]$,
by (\ref{EQ8175}) there is $\eta \in \Pc (\P )$, where $\eta :[a,q] \rightarrow (\cdot ,p]$
and the mapping $\p :D\rightarrow P$ defined by $\p =\f \cup \eta \upharpoonright (D\setminus \dom \f)$ is as desired.

{\it Case 3}: $\dom \f \cap (\cdot ,a )= \dom \f \cap (a,\cdot )=\emptyset$.
Then $a \parallel \dom \f$, and, hence, the set $D:=\dom \f \cup \{ a\}$ is convex.
By (\ref{EQ8174}) there is $b\in P\setminus [p,q] $
so the function $\p :=\f \cup \{ \la a,b\ra\}$ is an injection and $\dom \p =D$.
Since $a \parallel \dom \f$ we have $\p \in \Pc (\P )$.
Since $\P$ is upwards and downwards directed, there are $p' < a,b,p$ and $q' > a,b,q$,
which gives $\dom \p \cup \ran \p\subset (p',q')$.
Thus $\p \in \Pi$, $\f \subset \p$ and $a\in \dom \p$, so (bf1) is true.

(bf2) Let $\f \in \Pi$, $b\in P\setminus \ran \f$ and  $\dom \f \cup \ran \f \subset (p,q)$.
By (\ref{EQ8174}) there is $a\parallel [p,q] $; so, the set $D:=\dom \f \cup \{ a\}$ is convex.
The function $\p :=\f \cup \{ \la a,b\ra\}$ is an injection and $\dom \p =D$.
Since $a \parallel \dom \f$ we have $\p \in \Pc (\P )$.
Since $\P$ is upwards and downwards directed, there are $p' < a,b,p$ and $q' > a,b,q$,
which gives $\dom \p \cup \ran \p\subset (p',q')$.
Thus $\p \in \Pi$, $\f \subset \p$ and $b\in \ran \p$, so (bf2) is true.

By Theorem \ref{T8196} it remains to be proved that the poset $\la \Pi ,\supset\ra$ is $\k$-closed.
So let $\g \in \k$ and let $\f _\a \in \Pi$, where $\f _0 \subset \f _1 \subset \dots$.
Then, clearly,  $\f ^* \subset \f :=\bigcup _{\a <\g}\f _\a \in \Pc (\P)$
and the set $\dom \f :=\bigcup _{\a <\g}\dom \f _\a$ is convex.
By (\ref{EQ8170}) for each $\a <\g$ there are $p_\a,q_\a \in P$ such that $p_\a <\dom \f _\a \cup\ran \f _\a < q_\a$
and, since $\P$ is $<\k$-directed and $|\g|<\k$,
there are $p,q \in P$ such that $p< p_\a\leq q_\a < q$, for all $\a <\g$.
So we have $\dom \f  \cup\ran \f  =\bigcup _{\a <\g}\dom \f _\a \cup\ran \f _\a \subset (p,q)$;
thus $\f \in \Pi$ and $\f$ is a lower bound for the sequence $\la \f _\a :\a <\g\ra$.
\kdok
Generally speaking, the reversibility of a direct product $\prod _{i\in I}\X _i$ of $L$-structures
does not imply the reversibility of its factors $\X_i$, $i\in I$.
For example, if $\X =\la X,\r\ra$ is a non-reversible binary structure and $\Y =\la Y,\emptyset\ra$,
then the product $\X \times \Y =\la X\times Y, \emptyset\ra$ is a reversible structure and its factor $\X$ is not.
The following fact, used in the sequel, shows that the implication holds in some situations.
\begin{fac}\label{T8216}
If the direct product $\prod _{\a < \k }\X _\a$ of reflexive partial orders is reversible,
then all its factors $\X_\a$, $\a < \k$, are reversible.
\end{fac}
\dok
Let $\X _\a =\la X_\a ,\leq _\a\ra$, for $\a < \k$, and $\X :=\prod _{\a < \k }\X _\a =\la \prod _{\a < \k }X _\a ,\leq \ra$.
We prove the contrapositive of the statement.
Suppose that some of the factors, say $\X _0$, is not reversible;
let $f_0 \in \Cond (\X _0)\setminus \Aut (\X _0)$ and let $a,b\in X_0$, where $a\not\leq _0 b$ and $f_0(a)\leq _0 f_0(b)$.
Then the direct product of bijections $F:= f_0 \times \prod _{0<\a < \k }\id _{X_\a}: \prod _{\a < \k }X _\a \rightarrow \prod _{\a < \k }X _\a$ is a bijection
and $F(\la x_\a :\a <\k\ra)=\la f_0 (x_0), x_1, x_2, \dots \ra$.
If $x=\la x_\a :\a <\k\ra, y=\la y_\a :\a <\k\ra \in \prod _{\a < \k }X _\a$ and $x\leq y$,
then for each $\a <\k$ we have $x_\a \leq _\a y_\a$.
So $f_0 (x_0)\leq _0 f_0 (y_0)$ and, hence,  $\la f_0 (x_0), x_1, x_2, \dots \ra \leq \la f_0 (y_0), y_1, y_2, \dots \ra$,
that is $F(x)\leq F(y)$; thus $F\in \Cond (\X )$.

Let us choose $c_\a \in X_\a$, for $0<\a < \k$,
and let $x=\la a,c_1,c_2, \dots\ra$ and $y=\la b,c_1,c_2, \dots\ra$.
Since $a\not\leq _0 b$ we have $x\not\leq y$
but, since $f_0(a)\leq _0 f_0(b)$ we have $F(x)=\la f_0(a),c_1,c_2, \dots\ra \leq \la f_0(b),c_1,c_2, \dots\ra =F(y)$.
Thus $F\in \Cond (\X )\setminus \Aut (\X )$ and $\X$ is not reversible.
\hfill $\Box$
\begin{prop}\label{T8210}
The direct powers $\Q ^\k$ and $\Z ^\k$ are not reversible, whenever $\k \geq 2$.
\end{prop}
\dok
Since direct products are associative, by Fact \ref{T8216} it is sufficient to prove that the squares $\Q ^2$ and $\Z ^2$ are not reversible.

We use Theorem \ref{T8209} for $\k =\o$.
Both posets are lattices of size $\o$ and, thus they are upwards and downwards $<\o$-directed.
If $p=\la p_0,p_1\ra,q=\la q_0,q_1\ra\in \Q ^2$ and $p\leq q$, then $[p,q]$ is a rectangle in the rational plane,
$a=\la p_0-1,q_1+1\ra \parallel [p,q]$ and we have (\ref{EQ8174}).
If $r\in \Q ^2$, using the translation $f$ of the plane which moves $q$ to $r$
we obtain a partial condensation $\eta :=f\upharpoonright [p,q]:[p,q] \rightarrow (\cdot ,r]$;
using the translation $g$ of the plane which moves $p$ to $r$
we obtain a partial condensation $\eta :=g\upharpoonright [p,q]:[p,q] \rightarrow [r,\cdot )$
so (\ref{EQ8175}) and (\ref{EQ8176}) are true. For $\Z ^2$ we proceed in a similar way.
\kdok
We note that the same holds if instead of $\Q ^\k$ we regard direct powers of homogeneous universal linear orders (Hausdorff's $\eta _\a$-sets).
Proposition \ref{T8210} also shows that the direct product of two reversible structures (moreover, linear orders) can be non-reversible.

\paragraph{Acknowledgement.}
This research was supported by the Science Fund of the Republic of Serbia,
Program IDEAS, Grant No.\ 7750027:
{\it Set-theoretic, model-theoretic and Ramsey-theoretic
phenomena in mathematical structures: similarity and diversity}--SMART.

\end{document}